\makeatletter
\expandafter\def\csname ver@cleveref.sty\endcsname{2099/01/01 v9.9 stub}
\expandafter\def\csname opt@cleveref.sty\endcsname{capitalize,nameinlink}
\providecommand{\crefname}[3]{}
\providecommand{\Crefname}[3]{}

\providecommand{\crefformat}[2]{}
\providecommand{\Crefformat}[2]{}
\providecommand{\crefrangeformat}[2]{}
\providecommand{\Crefrangeformat}[2]{}
\providecommand{\crefmultiformat}[5]{}
\providecommand{\Crefmultiformat}[5]{}
\providecommand{\crefrangemultiformat}[5]{}
\providecommand{\Crefrangemultiformat}[5]{}
\providecommand{\crefdefaultlabelformat}[1]{}
\makeatother

\documentclass[final]{siamart171218}

\usepackage{graphicx}
\usepackage{multirow}
\usepackage{amsmath,amssymb,amsfonts}
\usepackage{mathtools}
\usepackage{mathrsfs}
\usepackage{xcolor}
\usepackage{booktabs}
\newtheorem{remark}[theorem]{Remark}

\headers{Correctness of Biot's model of in situ leaching}{A.~Meirmanov and A.~Senkebayeva}

\begin{document}

\title{Correctness of Biot's model of in situ leaching\\
for incompressible liquid and compressible solid components%
\thanks{Submitted to the editors.}}

\author{Anvarbek Meirmanov%
 \thanks{Institute of Physical and Technical Sciences,
    L.N.~Gumilyov Eurasian National University,
    Satbaev st., 2, Astana, Kazakhstan
    (\href{mailto:anvarbey1947@gmail.com}{anvarbey1947@gmail.com},
    \url{https://orcid.org/0000-0002-8543-3897}).
    }
  \and
 Akbota Senkebayeva%
  \thanks{School of Applied Mathematics,
    Kazakh-British Technical University,
    Tole bi st., 59, Almaty, Kazakhstan
    (\href{mailto:akbota.senkebayeva@gmail.com}{akbota.senkebayeva@gmail.com},
    \url{https://orcid.org/0000-0001-6646-5415}).}
}

\maketitle

\begin{abstract}
We study a mathematical model of in situ leaching of rare metals, in which
the joint filtration of two liquids is governed by the microscopic model
$\mathbb{A}^{\varepsilon}$. A key difficulty is the unknown
(\textbf{\emph{free}}) boundary $\Gamma(r)$ between solid and liquid components, determined by an additional condition on $\Gamma(r)$; no
standard methods exist for this nonlinear problem.
To resolve it we apply the \textbf{\emph{fixed point theorem}}.
For a given function $r(\boldsymbol{x},t)$ from a set
$\mathfrak{M}_{(0,T)}$ of sufficiently smooth functions describing the
skeleton structure, we consider the auxiliary problem
$\mathbb{B}^{\varepsilon}(r)$: an elliptic system for displacements of
the liquid and solid components coupled with parabolic equations for the
acid concentration. Selecting the \textbf{\emph{weak solution}} of minimal
smoothness, we apply the \textbf{\emph{homogenization method}} to pass from
the microscopic to the macroscopic description.
The resulting macroscopic model $\mathbb{H}(r)$ contains a homogenized
boundary condition that expresses the normal boundary velocity
$V_{N}=\partial r/\partial t$ as a linear function of the acid
concentration $c$. Since $c$ depends on $r$ via an operator
$\mathbb{F}\colon\mathfrak{M}_{(0,T)}\to\mathfrak{M}_{(0,T)}$,
we prove that $\mathbb{F}$ is Lipschitz continuous and, by Banach's
theorem, possesses a unique fixed point $r^{*}$, which yields the
unique solution $\mathbb{H}=\mathbb{H}(r^{*})$.
\end{abstract}



\section{Introduction}

In the proposed manuscript we consider a mathematical model of in situ
leaching of rare metals and prove existence and uniqueness theorems for the classical solvability of the initial-boundary value problem for the corresponding system of differential equations at the macroscopic level.

The extraction of rare metals by leaching is a very important task for our national economy. Natural deposits of uranium, nickel and other rare metals are complex geologically heterogeneous objects. Inhomogeneity means that the properties of ore bodies vary in space. Analyses of wells and cores show that the geological properties (porosity, permeability, etc.) of ore bodies are heterogeneous even within a single deposit. Very often the insufficient consideration of the consequences of inhomogeneities at the stage of operation planning becomes apparent too late, when the acid solution injected into the soil through injection wells appears far from the intended location. In addition, an important role is played by the concentration of the injected acid, the injection modes of acid solutions and other factors.

Hence, understanding the dynamics of fluids in heterogeneous porous media and the mechanism of dissolution of rocks by acids is of fundamental importance for the effective management of rare metal mining. This is achieved by creating a \textbf{\emph{Prototype of a hydrodynamic simulator}} for the ore body based on the corresponding mathematical models and a set of computer programs that visualise the physical process and determine the dynamic changes in the main characteristics of the mathematical model, which allows one to optimize the entire technological process. A \textbf{\emph{Hydrodynamic simulator of the ore body}} is a prototype of the hydrodynamic simulator supplemented with digital characteristics of the geometric and physical properties of the solid skeleton. The aim of this manuscript is to create a mathematical model for a prototype hydrodynamic simulator describing ore deposits with a given structure of the deposit's solid skeleton.

Currently, there is a wide range of mathematical models for describing the dynamics of in-situ leaching that can describe the physical processes under consideration, but only at the \textbf{\emph{Macroscopic Level}} (see~\cite{1,2,3} and references therein). Unlike \textbf{\emph{Microscopic}} models (the characteristic size is approximately tens of microns) in \textbf{\emph{Macroscopic}} models the characteristic size is decimeters or tens of decimeters. Because of this, the aforementioned models do not account for the microstructure of a continuous medium, since in such a model at each point the medium contains both a solid skeleton and liquid in the pores or cracks of this skeleton.

All these models are built on the same principle. Fluid dynamics is usually controlled by Darcy's system of filtration, or some modification of it. The equations describing the migration of acid and chemical reaction products are simply postulated and, roughly speaking, are some modifications of the diffusion-convection equations for the corresponding concentrations. The main thing in these postulates is the type of coefficients of the equations.

It is here that there is a wide variety of models, depending on the tastes and preferences of their authors. It is quite understandable, since the main mechanism of the physical process is focused on the unknown (\textbf{\emph{Free}}) boundary between the pore space and the solid skeleton and is not spelled out in any way in the proposed macroscopic models. This is where the rocks dissolve, changing the concentration of the injected acid, and this is where products appear inside the carrier liquid. Moreover, during the process, the geometry of the pore space (the geometry of the boundary separating the solid skeleton and the pore space) changes in time and space.

These fundamentally important changes occur at the microscopic level, corresponding to the average size of pores or cracks in rocks, while all of the proposed macroscopic models operate with completely different orders of scales and, therefore, do not distinguish the free boundaries, nor the features of the interaction of the acid with the rocks, which explains the wide variety of macroscopic mathematical models.

The authors of such models simply do not have an exact method for describing physical processes at the microscopic level based on the fundamental laws of classical mechanics of continuous media and chemistry, nor do they have the ability to take the microstructure into account in their macroscopic models. Therefore, they are forced to limit themselves to certain speculative considerations.

Due to this, a natural question arises. If there are several macroscopic models describing the same physical process under the same conditions, which of them most adequately reflects this process? Where is the criterion of adequacy here? It does not make sense to talk about an experiment, since each model has enough free parameters that are not related to the reservoir geometry (e.g., porosity), or to the physical characteristics of the process (such as the viscosity of the filtered liquids, or the physical properties of the solid skeleton and the like). So, with the variation of these free parameters, one may get a match with any experiment.

R.~Burridge and J.~B.~Keller~\cite{4} and E.~Sánchez-Palencia~\cite{5} were the first to explain that the exact description of filtration of liquids and seismic waves in rocks at the macroscopic level is possible if and only if:

(a) the physical process under consideration is described at the microscopic level by equations of Newtonian classical continuum mechanics (exact model);

(b) a set of small dimensionless parameters is selected.

(c) the macroscopic mathematical model is an exact asymptotic limit (homogenization) of exact mathematical models at the microscopic level, when the selected small parameters tend to zero.

Various special cases of exact macroscopic models of acoustics and fluid filtration in rocks  have been investigated by many authors (see~\cite{6,7,8,9,10,11,12,13}).  All these authors used different methods of homogenization, and the application of each of them was not an easy task requiring considerable effort and ingenuity.
Everything changed after the appearance of G.~Nguetseng's paper~\cite{14}, where the author proposed the \textbf{\emph{Method of two-scale convergence}} in periodic structures. What used to be  an art has become an ordinary routine, a matter of reference to the method. So, the homogenization theory has ceased to be an independent part of mathematical analysis (or the theory of differential equations) and the main efforts in homogenization have moved from theory to applications in mechanics, physics, biology, etc.

\section{Preliminaries}

\subsection{The problem statement}

As we have already noted, the derivation of macroscopic mathematical models should be based on the most exact mathematical model of the physical process at the microscopic level, described by the laws of Newtonian classical continuum mechanics~\cite{15} (see Appendix~A, section~A.7 in~\cite{16}).

Next, we will assume that the liquid in the pores is incompressible, and
the solid skeleton is a compressible medium.

The liquid motion in the pore space $\Omega^{\varepsilon}_{f}\subset\Omega$ for $t>0$ in the dimensionless variables~\cite{16} is governed by the stationary Stokes equations for an incompressible viscous fluid
\begin{equation}\label{eq2.1}
\nabla\cdot\mathbb{P}^{\varepsilon}_{f}=\nabla{p}^{0},\,\,
\mathbb{P}^{\varepsilon}_{f}=
\alpha^{\varepsilon}_{\mu}\mathbb{D}(x,\boldsymbol{v}^{\varepsilon}_{f})-
({p}^{\varepsilon}_{f}-{p}^{0})\mathbb{I},\,\,\,
\boldsymbol{v}^{\varepsilon}_{f}=\frac{\partial\boldsymbol{w}^{\varepsilon}_{f}}{\partial t},
\end{equation}
linearized continuity equation
\begin{equation}\label{eq2.2}
\nabla\cdot\boldsymbol{v}^{\varepsilon}_{f}=0
\end{equation}
for dynamic characteristics $\boldsymbol{w}^{\varepsilon}_{f}$ (\textbf{\emph{Liquid displacements}}), $\boldsymbol{v}^{\varepsilon}_{f}=
\frac{\partial\boldsymbol{w}^{\varepsilon}_{f}}{\partial t}$ (\textbf{\emph{Liquid velocity}}) and $p^{\varepsilon}_{f}$ (\textbf{\emph{Liquid pressure}}).

The motion of the compressible solid skeleton in the domain $\Omega^{\varepsilon}_{s}$ for $t>0$ is described by the stationary Lam\'{e} equations
\begin{equation}\label{eq2.3}
\nabla\cdot\mathbb{P}^{\varepsilon}_{s}=\nabla{p}^{0},\,\,
\mathbb{P}^{\varepsilon}_{s}=
\lambda_{0}\mathbb{D}(x,\boldsymbol{w}^{\varepsilon}_{s})-({p}^{\varepsilon}_{s}-{p}^{0})\mathbb{I},\,\,\,
\boldsymbol{v}^{\varepsilon}_{s}=\frac{\partial\boldsymbol{w}^{\varepsilon}_{s}}{\partial t},
\end{equation}
linearized continuity equation
\begin{equation}\label{eq2.4}
\frac{1}{c^{2}_{s}}({p}^{\varepsilon}_{s}-{p}^{0})+\nabla\cdot\boldsymbol{w}^{\varepsilon}_{s}=0
\end{equation}
for dynamic characteristic $\boldsymbol{w}^{\varepsilon}_{s}(\boldsymbol{x},t)$ (\textbf{\emph{Solid displacements}}) and $p^{\varepsilon}_{s}$ (\textbf{\emph{Solid pressure}}).

Diffusion of the acid and the products of chemical reactions in the pore space for $t>0$ is described by the diffusion equation
\begin{equation}\label{eq2.5}
\frac{\partial c^{\varepsilon}}{\partial t}=\nabla\cdot(\alpha_{0}\,\nabla\,c^{\varepsilon})
\end{equation}
for the acid \textbf{\emph{Acid concentration}} ${c}^{\varepsilon}$ and the transport equations
\begin{equation}\label{eq2.6}
\frac{\partial c^{\varepsilon}_{j}}{\partial t}+
\boldsymbol{v}^{\varepsilon}_{f}\cdot\nabla\,c^{\varepsilon}_{j}=0,\,\,j=1,...,k
\end{equation}
for \textbf{\emph{Concentrations of the products of chemical reactions}}
$c^{\varepsilon}_{j}(\boldsymbol{x},t),\,\,j=1,...,k$.

At the free boundary $\Gamma^{\,\varepsilon}$ between the liquid and solid components, the following boundary conditions
\begin{equation}\label{eq2.7}
\boldsymbol{w}^{\varepsilon}_{f}=\boldsymbol{w}^{\varepsilon}_{s},
\end{equation}
\begin{equation}\label{eq2.8}
\mathbb{P}^{\varepsilon}_{f}<\textbf{N}^{\varepsilon}>=
\mathbb{P}^{\varepsilon}_{s}<\textbf{N}^{\varepsilon}>,
\end{equation}
\begin{equation}\label{eq2.9}
(D^{\varepsilon}_{N}+\beta^{\varepsilon})\,c^{\varepsilon}+
\alpha_{0}\,\frac{\partial c^{\varepsilon}}{\partial\,N}=0,
\end{equation}
\begin{equation}\label{eq2.10}
(D^{\varepsilon}_{N}-v^{\varepsilon}_{f,N})\,c^{\varepsilon}_{j}=0,\,\,j=1,...,k,
\end{equation}
expressing the laws of conservation of mass and momentum~\cite{15} (see Appendix~A, section~A.7 in~\cite{16}), are satisfied.

Finally, an additional boundary condition is introduced, which should allow us to find a free
boundary that reduces the diameter of the solid components. We use the postulate in~\cite{16}, where the growth of biological tissue particles in a nutrient medium is considered, but with the opposite sign:
\begin{equation}\label{eq2.11}
D^{\varepsilon}_{N}(\boldsymbol{x},t)=\alpha^{\varepsilon}\,{c}^{\varepsilon}(\boldsymbol{x},t),\,\,\,
\boldsymbol{x}\in\,\Gamma^{\varepsilon},\,\,t>0.
\end{equation}
At the given boundaries with injection wells ${S}^{1}$ and producing wells ${S}^{2}$, and at the impermeable boundary ${S}^{0}$, the following conditions
\begin{equation}\label{eq2.12}
\mathbb{P}^{\varepsilon}<\boldsymbol{n}>=0,\,\,\boldsymbol{x}\in\,{S}^{1}\cup{S}^{2},\,\,t>0,
\end{equation}
\begin{equation}\label{eq2.13}
\frac{\partial}{\partial n}({c}^{\varepsilon}-{c}^{0})(\boldsymbol{x},t)=0,\,\,
\boldsymbol{x}\in S^{0},\,\,t>0,
\end{equation}
\begin{equation}\label{eq2.14}
\boldsymbol{w}^{\varepsilon}_{f}(\boldsymbol{x},t)=\boldsymbol{w}^{\varepsilon}_{s}(\boldsymbol{x},t)=0,
\,\,\,\boldsymbol{x}\in S^{0},\,\,t>0,
\end{equation}
\begin{equation}\label{eq2.15}
{c}(\boldsymbol{x},t)=c_{0}(\boldsymbol{x}),
\,\,\,\boldsymbol{x}\in S^{1}\cup S^{2},\,\,t>0,
\end{equation}
\begin{equation}\label{eq2.16}
{c}_{j}(\boldsymbol{x},t)=0,\,j=1,...,k,\,\,\,
\boldsymbol{x}\in S^{1},\,\,t>0
\end{equation}
are met.

The problem is completed with initial conditions
\begin{equation}\label{eq2.17}
{c}(\boldsymbol{x},0)=c_{0}(\boldsymbol{x}),\,\,\,
\boldsymbol{x}\in\Omega^{0}_{f},
\end{equation}
\begin{equation}\label{eq2.18}
\Gamma^{\varepsilon}_{t=0}=\Gamma^{0},\,\,\,
r(\boldsymbol{x},0)=r_{0}(\boldsymbol{x}),\,\,\boldsymbol{x}\in\Omega,
\end{equation}
\begin{equation}\label{eq2.19}
{c}_{j}(\boldsymbol{x},0)=0,\,j=1,...,k,\,\,\,
\boldsymbol{x}\in\Omega^{0}_{f}.
\end{equation}
In (\ref{eq2.1}) -- (\ref{eq2.19}) $\varepsilon$ is a small parameter equal to $\frac{L}{n}$,
where $n$ is an integer number, $n\gg{L}$, $D^{\varepsilon}_{N}$ is the \textbf{\emph{Normal velocity}} of the boundary $\Gamma^{\varepsilon}$ in the direction of the unit normal $\boldsymbol{N}^{\varepsilon}$ to the boundary $\Gamma^{\varepsilon}$, outward to the liquid domain $\Omega^{\varepsilon}_{f}$,
$\mathbb{D}(x,\boldsymbol{v}^{\varepsilon})=\frac{1}{2}\big(\nabla\boldsymbol{v}^{\varepsilon}+
(\nabla\boldsymbol{v}^{\varepsilon})^{*}\big)$ is a symmetric gradient of $\nabla\boldsymbol{v}^{\varepsilon}$,
$\mathbb{P}^{\varepsilon}=\chi^{\varepsilon}\mathbb{P}^{\varepsilon}_{f}+
(1-\chi^{\,\varepsilon})\mathbb{P}^{\varepsilon}_{s}$
is a \textbf{\emph{Stress tensor}}, $\mathbb{P}^{\varepsilon}_{f}$ is a stress tensor in the liquid component, $\mathbb{P}^{\varepsilon}_{s}$ is a stress tensor in the solid component, $p^{0}(\boldsymbol{x})$ is a given pressure, $p^{\varepsilon}_{f}$ and $p^{\varepsilon}_{s}$ are the pressures in the liquid and solid components, respectively. $p^{0}(\boldsymbol{x})=p^{i}$=\emph{const} for $\boldsymbol{x}\in{S}^{1}\cup{S}^{2}$, ${v}^{\varepsilon}_{f,N}=\chi^{\varepsilon}(\boldsymbol{v}^{\varepsilon}_{f}
\cdot\boldsymbol{N}^{\varepsilon})$ is the normal component of the liquid velocity at the free boundary, $c_{s}$ is the speed of sound in the solid skeleton 
and $\chi^{\varepsilon}$ is the characteristic function of the pore space. 

The positive constants $p^{1}$, $p^{2}$, $c_{s}$, $\alpha_{0}$ and $\lambda_{0}$ are supposed to be given.

The absolutely rigid solid skeleton was considered in~\cite{18}, where the key role in the microscopic description was played by the function $r(\boldsymbol{x},t)$ from the set 
\begin{multline}\label{eq2.20}
\mathfrak{M}_{(0,T)}=\{r\,
\in \mathbb{H}^{2+\gamma,\frac{2+\gamma}{2}}(\overline{\Omega}_{T}),\,\,
0<r(\boldsymbol{x})<\frac{1}{2},\,\,\,
-\theta\,\leqslant\,\frac{\partial{r}}{\partial{t}}
(\boldsymbol{x},t)\leqslant\,0,
\\
0<\gamma<1,\,\,\,\theta=\mbox{const}>0;\,\,\,
|r|^{(2+\gamma)}_{\Omega_{T}}\leqslant\,M_{0}\},
\end{multline}
which determined the structure of the pore space.

We call the problem (\ref{eq2.1}) -- (\ref{eq2.5}), (\ref{eq2.7}) -- (\ref{eq2.9}), (\ref{eq2.11}) -- (\ref{eq2.15}), (\ref{eq2.17}), (\ref{eq2.18}) a problem $\mathbb{A}^{\varepsilon}$.

As the problem $\mathbb{B}^{\varepsilon}(r)$ we denote the problem $\mathbb{A}^{\varepsilon}$ without the initial condition for the free boundary and the additional boundary condition (\ref{eq2.11}) at the free boundary, which defines this boundary, but with a known structure of the pore space, given by the function $r(\boldsymbol{x},t)$ from the set $\mathfrak{M}_{(0,T)}$ coupled with an additional term in the dynamic equation for the liquid component in the form
\begin{equation}\label{eq2.21}
\nabla\cdot\mathbb{P}^{\varepsilon}_{f}=
\varepsilon^{2}\frac{\partial\boldsymbol{w}^{\varepsilon}_{f}}{\partial{t}}+\nabla{p}^{0}.
\end{equation}
As we have mentioned in the abstract, we first consider the problem $\mathbb{B}^{\varepsilon}(r)$ with a given structure of the pore space. In this problem, for a fixed $\varepsilon>0$, the solid skeleton is a union of disjoint sets sufficiently close to balls of radius $\varepsilon\,r$, slowly decreasing in volume, which simplifies the geometry of the original pore space and allows us to prove the existence of approximate solutions. As usual, almost every new problem has multiple choices. For example, for our case we may consider non-stationary Stokes equations, but then we must somehow find a priori estimates for the liquid velocities keeping in mind the difficulties with the free boundary, separating the liquid and solid components.

We refer to the problem (\ref{eq2.1})--(\ref{eq2.4}), (\ref{eq2.7}), (\ref{eq2.8}), (\ref{eq2.12}), (\ref{eq2.14}), (\ref{eq2.18}) as the \textbf{\emph{dynamic problem}} $\mathbb{B}^{\varepsilon}(r)$,  problem (\ref{eq2.5}), (\ref{eq2.9}), (\ref{eq2.13}), (\ref{eq2.15}), (\ref{eq2.17}) as
\textbf{\emph{Diffusion problem}} $\mathbb{B}^{\varepsilon}(r)$, the homogenization of the dynamic problem  $\mathbb{B}^{\varepsilon}(r)$ as a homogenized dynamic problem $\mathbb{H}(r)$, and homogenization of the diffusion problem $\mathbb{B}^{\varepsilon}(r)$ as homogenized diffusion problem $\mathbb{H}(r)$.

Homogenization of the problem $\mathbb{B}^{\varepsilon}(r)$ gives us a homogenized model $\mathbb{H}(r)$, consisting of the homogenized dynamic problem $\mathbb{H}(r)$ and the homogenized diffusion problem $\mathbb{B}^{\varepsilon}(r)$.

Finally, homogenizing the boundary condition (\ref{eq2.11}) yields an operator $\mathbb{F}$ that transforms the set $\mathfrak{M}_{(0,T)}$ into itself. The unique fixed point $r^{*}$ of this operator
will determine the desired unique homogenization $\mathbb{H}=\mathbb{H}(r^{*})$ of problem $\mathbb{A}^{\varepsilon}$.

To homogenize the dynamic problem $\mathbb{B}^{\varepsilon}(r)$ for the liquid component, the dynamic problem $\mathbb{B}^{\varepsilon}(r)$ for the solid component, and the diffusion problem $\mathbb{B}^{\varepsilon}(r)$, we will use a modification of Nguetseng's two-scale convergence method adapted for \textbf{\emph{Structures with special periodicity}}~\cite{18}.

The physical process we are considering is quite long-lasting (the fluid filtration rate is several meters per year). Therefore, theorems of greatest interest are those on the existence of solutions to the corresponding initial-boundary value problems globally in time. On the other hand, due to the strong nonlinearity of free boundary problems (\cite{18,19}), it is usually impossible to prove any results globally in time for mathematical models at the microscopic level. Thus, the only possible results are theorems on the existence of a generalized or classical solution locally in time.

It should be noted that the Stefan problem was formulated only at the macroscopic level, and for a one-phase problem where the free boundary is monotone in time, it allowed one to prove the existence of a classical solution globally in time~\cite{20}.

But in the general situation, how can we obtain macroscopic mathematical models if we know nothing about the existence of solutions of the microscopic mathematical models of which they are supposed to be the limit?

To get around these difficulties, we follow the ideas of~\cite{18} and will use the fixed point theorems~\cite{21}. To do this, we define the structure of the pore space, given by the characteristic function 
\begin{equation*}
\chi^{\varepsilon}(\boldsymbol{x},t)=\chi\big(r(\boldsymbol{x},t);\boldsymbol{y}\big)=
\chi\big(r(\boldsymbol{x},t);[\frac{\boldsymbol{x}}{\varepsilon}]\big),
\end{equation*}
periodic in the variable $\boldsymbol{y}$. 

As we have already noted, in the general case, solving the emerging problem is almost impossible. Therefore, it is reasonable to limit ourselves to the simplest cases. For example, when a non-negative function $r(\boldsymbol{x},t)$ from the set $\mathfrak{M}_{(0,T)}$ uniquely determines the characteristic function of the pore space $\chi(r;\boldsymbol{y})$.

That is, for a fixed $r\in \mathfrak{M}_{(0,T)}$ we consider the initial boundary value problem for the model $\mathbb{B}^{\varepsilon}(r)$ in a given domain $\Omega^{\varepsilon}_{f}(r)$ occupied by the liquid component and a given domain $\Omega^{\varepsilon}_{s}(r)$ occupied by the solid component for determining the main characteristics of the medium (velocities, displacements, pressures, and acid concentration), which together amounts to the initial boundary value problem $\mathbb{A}^{\varepsilon}$ without the additional boundary condition (\ref{eq2.11}). In this problem for a fixed $\varepsilon>0$ the solid skeleton is a union of disjoint sets sufficiently close to balls of radius $\varepsilon{r}$, slowly decreasing in volume, which simplifies the geometry of the original pore space and allows us to prove the existence of approximate solutions.

To understand what the \textbf{Homogenized Problem} $\mathbb{H}(r)$  of the problems $\mathbb{B}^{\varepsilon}(r)$ should be, a formal homogenization of the problem $\mathbb{A}^{\varepsilon}$ is performed beforehand. The conditions for the existence of homogenization of the boundary condition (\ref{eq2.11}) are formulated in Lemma~\ref{l2.8}. 

If $r^{\varepsilon}(\boldsymbol{x},t)$ defines the structures of the solid skeleton and pore space in the problem $\mathbb{A}^{\varepsilon}$ and $r^{\varepsilon}\rightarrow r^{*}$ as $\varepsilon\rightarrow 0$, then the homogenized problem $\mathbb{H}(r^{*})$ of the problem $\mathbb{B}^{\varepsilon}(r^{*})=\mathbb{A}^{\varepsilon}$ should coincide with the homogenization $\mathbb{H}$ of the problem $\mathbb{A}^{\varepsilon}$ without the homogenization of the boundary condition (\ref{eq2.11}).

It is clear that the homogenization of the boundary condition (\ref{eq2.11}) with a given structure of the pore space, defined by the function $r\in\mathfrak{M}_{(0,T)}$, forms an operator $\mathbb{F}:\mathfrak{M}_{(0,T)}\rightarrow\mathfrak{M}_{(0,T)}$, whose unique fixed point $r^{*}(\boldsymbol{x},t)$ determines the required unique homogenization $\mathbb{H}$ of the problem $\mathbb{A}^{\varepsilon}$.

To solve the problem $\mathbb{H}(r)$ we, first of all, have to solve the linear problem $\mathbb{B}^{\varepsilon}(r)$ and then find its homogenization $\mathbb{H}(r)$ as
$\varepsilon\rightarrow 0$.

In turn, the linear problem $\mathbb{B}^{\varepsilon}(r)$ is decomposed into a sequential solution of the \textbf{dynamic problem} $\mathbb{B}^{\varepsilon}(r)$, defining the dynamic characteristics
$\boldsymbol{w}^{\varepsilon}_{f}$, $\boldsymbol{v}^{\varepsilon}_{f}$, $\boldsymbol{w}^{\varepsilon}_{s}$,
$p^{\varepsilon}_{f}$, $p^{\varepsilon}_{s}$ and the \textbf{Diffusion Problem} $\mathbb{B}^{\varepsilon}(r)$ defining the acid concentration $c^{\,\varepsilon}$.

Due to the linearity of these problems, the existence and uniqueness of a weak solution to each of them follows from the corresponding a priori estimates and known methods for solving linear differential equations. For example, \textbf{\emph{Galerkin method}}~\cite{22}.

The next step is the homogenization of the resulting mathematical model $\mathbb{B}^{\varepsilon}(r)$. To homogenize this problem, we use the modified Nguetseng's two-scale convergence method~\cite{18}, which makes it quite easy to get a macroscopic mathematical model $\mathbb{H}(r)$. It consists of a dynamic model $\mathbb{H}(r)$ for the liquid component, a dynamic model $\mathbb{H}(r)$ for the solid component, and a diffusion model $\mathbb{H}(r)$ for the acid concentration. But since this method was developed only for the homogenization of functionals, we needed to write down the original mathematical model in the form of a system of integral identities equivalent to the original system of differential equations, supplemented by the corresponding boundary conditions.

The integral identities equivalent to the dynamic Stokes and Lam\'{e} equations, as well as the integral identities for the diffusion equation with standard boundary conditions are well known.

But the equivalent expression of differential equations in the form of integral identities is a general and rather difficult challenge for the free boundary problems. The happy exception was the Stefan 
problem (\cite{23,24}) describing phase transitions in pure (without impurities) media. Such, for example, as ``water-ice'' or chemically pure metals. The authors of the cited works O. A. Oleinik and S. L. Kamenomostskaya (S. Kamin) managed to reformulate the problem so that the equivalent formulation in the form of an integral identity, in the case of the existence of a classical solution, contained both the equation of thermal conductivity outside the free boundary and the condition itself on the free boundary. This approach allowed them to prove quite simply the existence and uniqueness of the weak solution to the integral identity under minimal conditions for the smoothness of the solution. At the same time, in the case of the existence of a classical solution of the Stefan problem, the latter must coincide with the weak solution. The question of the existence of a classical solution globally in time to the one-phase Stefan problem remained open until 1975~\cite{20}. The existence of a classical solution to the two-phase Stefan problem locally in time was proved in 1979~\cite{19}, provided that the modulus of the temperature gradient at the free boundary at the initial time is positive. Moreover, in~\cite{25} it was shown that if this condition is violated, the classical solution of the two-phase Stefan problem does not exist.

Given the above, in the mathematical model of in situ leaching that we proposed, it was very important to find an equivalent formulation of the problem in the form of a system of integral identities that would require minimal smoothness of the solutions to the problem, which was done.

But, for example, at the same time, for the acid concentration ${c}^{\varepsilon}$, defined only in pore space $\Omega^{\varepsilon}_{f,T}(r)=\bigcup_{t=0}^{t=T}\Omega^{\varepsilon}_{f}(r)$, it was necessary to find an extension of these solutions from the domain of definition onto the domain $\Omega_{T}=\Omega\times(0,T)$ while preserving their best differential properties. To do this, we used the results on the extension of functions, formulated in~\cite{26,27}, so that the extensions $c^{\,\varepsilon}(\boldsymbol{x},t)$ satisfy the integral identities of conservation of acid concentration, equivalent to the corresponding differential equations together with boundary and initial conditions.

A priori estimates of weak solutions (i.e., solutions of corresponding integral identities) usually require a special choice of test functions in the integral identity and integration by parts.  For the latter, sufficient smoothness of the boundary of the pore space $\Omega^{\varepsilon}_{f}(r)$ (the domain filled with fluid) is necessary. The smoothness of the boundary $\partial\Omega^{\varepsilon}_{f}(r)$ is determined by the function $r\in\mathfrak{M}_{(0,T)}$. This simple fact is central to the derivation of a priori estimates.

We will show that the operator $\mathbb{F}=\mathbb{F}(r)$ is Lipschitz continuous with the corresponding constant bounded by some linear function of $T$. This property allows us to prove the existence of a unique fixed point $r^{*}(\boldsymbol{x},t)$, at least locally in time.

Finally, using the smoothness of the solutions to the problem $\mathbb{H}(r)$ we prove the correctness (existence and uniqueness) of the mathematical model $\mathbb{H}$ for any $T>0$.

In our manuscript, we use the notation adopted in~\cite{28,29,30}.
\begin{remark}\label{Rm2.1}
It is clear that we may find concentrations of products of chemical reactions after finding the solutions to the problems $\mathbb{B}^{\varepsilon}(r)$ and $\mathbb{H}(r)$.
\end{remark}

\subsection{Dimensionless parameters}

The dimensionless parameter  $\varepsilon=\frac{l}{L}$ is taken as a small parameter.  Furthermore, without loss of generality, we assume that $\displaystyle\,\varepsilon=\frac{1}{n}$,
 where $n$ is an integer.

The dimensionless parameter $\alpha^{\varepsilon}_{\mu}$ characterizes the viscosity of the liquid in pores:
\begin{equation*}
\alpha^{\varepsilon}_{\mu}=\frac{2\,\mu}{L\,g\,\tau\,\rho_{\,0}}
\end{equation*}
and dimensionless parameter $\alpha^{\varepsilon}_{0}$ characterizes the speed of dissolution of the solid skeleton.

Diffusion of acid is characterized by dimensionless coefficient
\begin{equation*}
\alpha_{0}=\frac{D\,T}{L^2}.
\end{equation*}
Here $l$ is the characteristic pore size and $L$ is the characteristic size of the physical domain under consideration, $\tau$ is the characteristic duration time of the physical process, $\rho_{\,0}$ is the density of water, $g$ is the acceleration of gravity and $\mu$ is the dynamic viscosity of the liquid, $\varrho_{s}$ is the dimensionless density of the solid skeleton, related to the density of water $\rho_{0}$ and  $\varrho_{f}$ is the dimensionless density of the liquid component related to the density of water $\rho_{0}$ and $D$ is the acid diffusion coefficient. Parameters $\alpha^{\varepsilon}_{c}$ and $\alpha^{\varepsilon}_{\mu}$ may depend on the small parameter $\varepsilon$ and parameters $\alpha_{0}$, $\beta$ and $\beta_{j}\,\,j=0,1,...,k,$ are given positive constants that do not depend on the small parameter $\varepsilon$.

In the present manuscript we consider \textbf{\emph{Biot's model}} for in situ leaching, when
\begin{equation*}
\alpha^{\varepsilon}_{\mu}=\varepsilon^{2}\mu_{1},\,\,\mbox{and}\,\,\mu_{1}=\mbox{const}>0.
\end{equation*}

\subsection{The structure of the pore space}

In what follows all functions of the type $\varphi(\boldsymbol{y};\boldsymbol{x},t)$, where $(\boldsymbol{x},t)\in \Omega$ and $\boldsymbol{y}\in \mathbb{R}^{3}$
are considered 1-periodic in variable $\boldsymbol{y}$:
\begin{equation}\label{eq.2.22}
\varphi(\boldsymbol{y};\boldsymbol{x},t)=
\varphi\big(\boldsymbol{\varsigma}(\boldsymbol{y});\boldsymbol{x},t)\big),
\,\,\,
\boldsymbol{y}=[|\boldsymbol{y}|]+\varepsilon\,\boldsymbol{\varsigma}(\boldsymbol{y}),\,\,
[|\boldsymbol{y}|]=([|y_{1}|],\,[|y_{2}|],\,[|y_{3}|]).
\end{equation}
The number $[|a|]$ is the integer part of the number $a$.

For the problem $\mathbb{A}^{\varepsilon}$ we restrict ourselves with the
simplest structure of the pore space $\textbf{Y}$, where
\begin{multline}\label{eq2.23}
\textbf{Y}=\{\boldsymbol{y}\in\mathbb{R}^{3}:-\frac{1}{2}<y_{k}<\frac{1}{2},\,\,k=1,2,3\},
\\
\textbf{Y}_{s}(r^{*})=\{\boldsymbol{y}\in\textbf{Y}: |\boldsymbol{y}|=(y_{1}^{2}+y_{2}^{2}+y_{3}^{2})^{\frac{1}{2}}<r^{*}\},\,\,
\textbf{Y}_{f}(r^{*})=\{\boldsymbol{y}\in\textbf{Y}: |\boldsymbol{y}|>r^{*}\},
\\
\gamma(r^{*})=\partial\textbf{Y}_{f}(r^{*})\cap\partial\textbf{Y}_{s}(r^{*})
\end{multline}
and   
\begin{equation}\label{eq2.24}
\chi(r^{*};\boldsymbol{y})=\frac{\mbox{sgn}(|\boldsymbol{y}|-r^{*})+1}{2},\,\,
\chi^{\varepsilon}(\boldsymbol{x},t)=\chi\big(r(\boldsymbol{x},t);[\frac{\boldsymbol{x}}{\varepsilon}]\big).
\end{equation}
As $\displaystyle\,\boldsymbol{N}(r^{*})=-\frac{\boldsymbol{y}}{|\boldsymbol{y}|}$ we denote the outward unit normal to the domain ${Y}_{f}(r^{*})\subset{Y}$.

The same notations as for problems $\mathbb{A}^{\varepsilon}$ we use for the problems $\mathbb{B}^{\varepsilon}(r)$, where instead of $r^{*}$ we insert $r$:
\begin{multline}\label{eq2.25}
\textbf{Y}=\{\boldsymbol{y}\in\mathbb{R}^{3}:-\frac{1}{2}<y_{k}<\frac{1}{2},\,\,k=1,2,3\},
\\
\textbf{Y}_{s}(r)=\{\boldsymbol{y}\in\textbf{Y}: |\boldsymbol{y}|<r\},\,\,
\textbf{Y}_{f}(r)=\{\boldsymbol{y}\in\textbf{Y}: |\boldsymbol{y}|>r\},
\\
\gamma(r)=\partial\textbf{Y}_{f}(r)\cap\partial\textbf{Y}_{s}(r)
\\
\chi(r;\boldsymbol{y})=\frac{\mbox{sgn}(|\boldsymbol{y}|-r)+1}{2},\,\,
\chi^{\varepsilon}(\boldsymbol{x},t)=\chi\big(r(\boldsymbol{x},t);
[\frac{\boldsymbol{x}}{\varepsilon}]\big).
\end{multline}
In (\ref{eq2.25}) $\chi^{\varepsilon}(\boldsymbol{x},t)$ is a 
\textbf{\emph{Characteristic function of the solid skeleton}}. That is, $\chi^{\varepsilon}(\boldsymbol{x},t)=0$ in the solid skeleton and $\chi^{\varepsilon}(\boldsymbol{x},t)=1$ 
in the pore space.

\subsection{Domains and boundaries}

Let $Q\subset\mathbb{R}^{3}$. As $\overline{Q}$ we denote the closure of the domain $Q$.

Let also $\Omega\subset \mathbb{R}^{3}$ be a bounded domain with piecewise smooth boundary
$S=\partial\Omega=\bar{S}^{0}\cup\bar{S}^{1}\cup\bar{S}^{2}$.

The boundary $S^{0}\subset \mathbb{R}^{3}$ is impermeable to liquid in the pore space,
the boundary $S^{1}\subset \mathbb{R}^{3}$ simulates injection wells and the boundary
$S^{2}\subset \mathbb{R}^{3}$ simulates production wells.

In what follows we will assume that $\Omega$ is the unit cube, $\Omega_{T}=\Omega\times(0,T)\subset\mathbb{R}^{3}$,

$S^{0}=\{\boldsymbol{x}:x_{3}=\pm\frac{1}{2},\,-\frac{1}{2}\,\leqslant\,x_{1},\,x_{2}\,
\leqslant\,\frac{1}{2}\}$,

$S^{1}=\{\boldsymbol{x}:x_{1}=
-\frac{1}{2},\,-\frac{1}{2}\,\leqslant\,x_{2},\,x_{3}\,\leqslant\,\frac{1}{2}\}$,

$S^{2}=\{\boldsymbol{x}:x_{1}=\frac{1}{2},\,-\frac{1}{2}\,
\leqslant\,x_{2},\,x_{3}\,\leqslant\,\frac{1}{2}\}$,

Next we put

$\Omega^{\varepsilon}_{f}(r)=\{\boldsymbol{x}\in\Omega: \chi^{\varepsilon}(\boldsymbol{x},t)=1\}$,\,\,\,
$\Omega^{\varepsilon}_{s}(r)=\{\boldsymbol{x}\in\Omega: \chi^{\varepsilon}(\boldsymbol{x},t)=0\}$,
$\Omega^{\varepsilon}_{f,T}(r)=\bigcup_{t=0}^{t=T}\Omega^{\varepsilon}_{f}(r)$,\,\,\,
$\Omega^{\varepsilon}_{s,T}(r)=\bigcup_{t=0}^{t=T}\Omega^{\varepsilon}_{s}(r)$,\,\,\,
$\Gamma^{\varepsilon}(r)=\overline{\Omega}_{f}(r)\cap\overline{\Omega}_{s}(r)$,\,\,
$\Gamma^{\varepsilon}_{T}=\bigcup_{t=0}^{t=T}\Gamma^{\,\varepsilon}(r)$  and
\begin{multline*}
\Omega=\bigcup_{\boldsymbol{k}\in \mathbb{Z}}\overline{\Omega}^{\,\boldsymbol{k},\varepsilon},\,\,\,
\Omega^{\boldsymbol{k},\varepsilon}=\{\boldsymbol{x}\in\Omega:\boldsymbol{x}=
\varepsilon\boldsymbol{k}+\varepsilon\boldsymbol{y}\},\,\,
\Omega_{f}^{\boldsymbol{k},\varepsilon}(r)=\Omega^{\varepsilon}_{f}(r)\cap\Omega^{\,\boldsymbol{k},\varepsilon},
\\
\Omega_{s}^{\boldsymbol{k},\varepsilon}(r)=
\Omega^{\varepsilon}_{s}(r)\cap\Omega^{\,\boldsymbol{k},\varepsilon},\,
\Omega^{0}_{f}=\Gamma^{\boldsymbol{k},\varepsilon}(r)=
\Gamma^{\varepsilon}(r)\cap\Omega^{\,\boldsymbol{k},\varepsilon}
\end{multline*}
for all $\boldsymbol{k}=(k_{1},k_{2},k_{3})$, $k_{1},k_{2},k_{3}\in \mathbb{Z}$ (integer numbers) and

for all $\boldsymbol{y}\in\,Y=(-\frac{1}{2},\frac{1}{2})^{3}\subset\mathbb{R}^{3}$.

Thus,
\begin{multline*}
\Omega^{\varepsilon}_{f}(r)=\{\boldsymbol{x}\in\Omega:
\chi^{\varepsilon}(\boldsymbol{x},t)=1\},\,\,\Omega^{0}_{f}=\Omega^{\varepsilon}_{f}(r_{0}),
\\
\Omega^{\varepsilon}_{s}(r)=\{\boldsymbol{x}\in\Omega:
\chi^{\varepsilon}(\boldsymbol{x},t)=0\},\,\,\Omega^{0}_{s}=\Omega^{\varepsilon}_{s}(r_{0}),
\\
\Omega^{\boldsymbol{k},\varepsilon}_{j}(r)=
\Omega^{\boldsymbol{k},\varepsilon}\cap\Omega^{\varepsilon}_{j}(r),\,\,\,j=f,s,
\\
\Gamma^{\varepsilon}(r)=\overline{\Omega}^{\varepsilon}_{f}(r)
\bigcap\overline{\Omega}^{\varepsilon}_{s}(r)=\bigcup_{k=1}^{n^{3}}\Gamma^{\,\varepsilon,k}(r),\,\,
\Gamma^{\varepsilon,k}(r)=\Omega^{\boldsymbol{k},\varepsilon}\cap\Gamma^{\,\varepsilon}(r).
\end{multline*}
We call the structure, defined by the formula (\ref{eq2.25}) 
\textbf{\emph{Structure with special periodicity}}.

For any continuous in $\Omega_{f}$ function $u(\boldsymbol{x})$ its limits at the points $\boldsymbol{x}_{0}$ from $\Gamma^{\,\varepsilon}(r)$ are denoted as
\begin{multline*}
\,\,\,\,\,\,\,\,\,u(\boldsymbol{x}_{0}+0)=
\lim_{\boldsymbol{x}\rightarrow \boldsymbol{x}_{0}}u(\boldsymbol{x}),
\,\,\boldsymbol{x}\in \Omega^{\varepsilon}_{f,t_{0}}(r),\,\,\,
\boldsymbol{x}_{0}\in \Gamma^{\varepsilon}(r),\\
u(\boldsymbol{x}_{0}-0)=
\lim_{\boldsymbol{x}\rightarrow\boldsymbol{x}_{0}}
u(\boldsymbol{x}),\,\,\boldsymbol{x}\in\Omega^{\varepsilon}_{s,t_{0}}(r),\,\,\,
\boldsymbol{x}_{0}\in \Gamma^{\varepsilon}(r).
\end{multline*}
Next, we choose a small parameter $\displaystyle\,\varepsilon=\frac{1}{n},\,\, n=1,2,3,...$ so that the boundary condition (\ref{eq2.12}) on the boundary ${S}^{1}\cup{S}^{2}$ makes sense.

\subsection{Matrices, tensors and differential operators}

We fix the standard Cartesian orthogonal basis
$\boldsymbol{e}^{1},\,\boldsymbol{e}^{2},\,\boldsymbol{e}^{3}$
in $\mathbb{R}^{3}$ and $\mathbb{A}$, $\mathbb{B}$ and $\mathbb{C}$ are \textbf{\emph{Tensors}} (linear transformations $\mathbb{R}^{3}\rightarrow\mathbb{R}^{3}$).
The action of the tensor $\mathbb{A}$ on the vector $\boldsymbol{b}$
is denoted as the vector $\boldsymbol{c}=\mathbb{A}<\boldsymbol{b}>$. As $(\boldsymbol{a}\cdot\boldsymbol{b})$ we denote the \textbf{\emph{Scalar product}} of vectors $\boldsymbol{a}$ and $\boldsymbol{b}$.
The product $\mathbb{C}=\mathbb{A}\cdot\mathbb{B}$ is a transformation
$\mathbb{A}: \mathbb{B}(\mathbb{R}^{3})\rightarrow\mathbb{R}^{3}$, where
$\mathbb{B}(\mathbb{R}^{3})=\{\boldsymbol{y}\in\mathbb{R}^{3}:\boldsymbol{y}=
\mathbb{B}(\boldsymbol{x}),\,\forall \boldsymbol{x}\in\mathbb{R}^{3}\}$.
$\mathbb{I}$ is a unit tensor: $\mathbb{I}\cdot\mathbb{A}=\mathbb{A}\cdot\mathbb{I}=\mathbb{A}$ for any tensor $\mathbb{A}$.

For any vectors $\boldsymbol{a}$, $\boldsymbol{b}$, $\boldsymbol{c}$ as $\boldsymbol{a}\otimes\boldsymbol{b}$ we denote the \textbf{\emph{Diad}} (second-rank tensor), where $(\boldsymbol{a}\otimes\boldsymbol{b})<\boldsymbol{c}>=
\boldsymbol{a}(\boldsymbol{b}\cdot\boldsymbol{c})$.

As $\mathbb{J}^{ij}$ we denote the second-rank tensor $\mathbb{J}^{ij}=\frac{1}{2}
(\boldsymbol{e}^{i}\otimes\boldsymbol{e}^{j}+\boldsymbol{e}^{j}\otimes\boldsymbol{e}^{i})$.

Then $\displaystyle\,\mathbb{A}=\sum_{i,j=1}^{3}{a}_{ij}\boldsymbol{e}^{i}\otimes\boldsymbol{e}^{j}$.

The second order tensor $\mathbb{A}$ is symmetric, if
$(\mathbb{A}<\boldsymbol{e}^{j}>\cdot\boldsymbol{e}^{i})=
(\mathbb{A}<\boldsymbol{e}^{i}>\cdot\boldsymbol{e}^{j})$.

As can be easily seen, any second-rank tensor is a linear mapping from $\mathbb{R}^{3}\rightarrow\mathbb{R}$.

Next, we determine the fourth-rank tensor $\mathfrak{N}$ as bilinear mapping from $\mathbb{R}^{3}\rightarrow\mathbb{R}$.

By $({A})$, $({B})$ and $({C})$ we denote the corresponding  to tensors $\mathbb{A}$, $\mathbb{B}$ and $\mathbb{C}$ matrices in the chosen Cartesian coordinate system:
\begin{equation*}
({A})=\left(
\begin{array}{cccc}
a_{11}&a_{12}&a_{13}\\
a_{21}&a_{22}&a_{23}\\
a_{31}&a_{32}&a_{33}
\end{array}
\right),\,\,({B})=
\left(
\begin{array}{cccc}
b_{11}&b_{12}&b_{13}\\
b_{21}&b_{22}&b_{23}\\
b_{31}&b_{32}&b_{33}
\end{array}
\right),\,\,({C})=
\left(
\begin{array}{cccc}
c_{11}&c_{12}&c_{13}\\
c_{21}&c_{22}&c_{23}\\
c_{31}&c_{32}&c_{33}
\end{array}
\right),
\end{equation*}
For matrices, the usual operations of sum $({A})+({B})$, multiplication by scalars $\alpha({B})$ and product $(A)\cdot(B)$ are defined.

Let $\boldsymbol{w}(\boldsymbol{x},t)=
\big({w}_{1}(\boldsymbol{x},t),{w}_{2}(\boldsymbol{x},t),{w}_{3}(\boldsymbol{x},t)\big)$ and
$\mathbb{D}(x,\boldsymbol{w})=
\frac{1}{2}\big(\nabla_{x}\boldsymbol{w}+(\nabla_{x}\boldsymbol{w})^{*}\big)$.

Then the second-rank symmetric tensor
\begin{equation*}
\begin{aligned}
\mathbb{D}(x,\boldsymbol{w})
&= \frac{1}{2} \sum_{i,j=1}^{3} \Bigl(
d_{ij}(x,\boldsymbol{w})\, \boldsymbol{e}^{i} \otimes \boldsymbol{e}^{j}
+ d_{ji}(x,\boldsymbol{w})\, \boldsymbol{e}^{j} \otimes \boldsymbol{e}^{i}
\Bigr), \\
d_{ij}(x,\boldsymbol{w})
&= \frac{\partial w_i}{\partial x_j}(x,\boldsymbol{w}),
\qquad i,j = 1,2,3.
\end{aligned}
\end{equation*}
is called the \textbf{\emph{Symmetric gradient}} of the vector $\boldsymbol{w}$.

We put as the definition  
\begin{multline}\label{eq2.26}
\mathbb{D}(x,\boldsymbol{w})<\boldsymbol{a}>\,\,\stackrel{df.}{=}\frac{1}{2}\sum_{i,j=1}^{3}
\big({d}_{ij}(x,\boldsymbol{w})(\boldsymbol{e}^{i}\otimes\boldsymbol{e}^{j})+{d}_{ji}(x,\boldsymbol{w})
(\boldsymbol{e}^{j}\otimes\boldsymbol{e}^{i})\big)<\boldsymbol{a}>=
\\
\frac{1}{2}\sum_{i=1,j}^{3}({d}_{ij}(x,\boldsymbol{w}){a}_{j}\boldsymbol{e}^{i}+
{d}_{ji}(x,\boldsymbol{w}){a}_{i}\boldsymbol{e}^{j})=
\sum_{i=1}^{3}\frac{1}{2}\sum_{i=1,j}^{3}{a}_{i}\boldsymbol{e}^{j},
\\
\mathbb{D}(x,\boldsymbol{w}):\mathbb{D}(x,\boldsymbol{v})=\frac{1}{2}\sum_{i=1,j}^{3}
{d}_{ij}(x,\boldsymbol{w}){d}_{ij}(x,\boldsymbol{v}).
\end{multline}
Then
\begin{equation}\label{eq2.27}
\mathbb{D}(x,\boldsymbol{w}):\mathbb{D}(x,\boldsymbol{w})=
\sum_{i,j=1}^{3}{d}_{ij}(x,\boldsymbol{w}){d}_{ij}(x,\boldsymbol{w}),\,\,
|\mathbb{D}(x,\boldsymbol{w})|^{2}=\sum_{i,j=1}^{3}|{d}_{ij}(x,\boldsymbol{w})|^{2}.
\end{equation}

\subsection{Poincar\'{e} inequality}

\begin{lemma}\label{l2.1}
Let~\cite{31} $\Omega\subset\mathbb{R}^{3}$ be a bounded domain with a piecewise smooth Lipschitz boundary.

Then for any function $\boldsymbol{w}^{\varepsilon}\in\stackrel{\!\!\circ}{\mathbb{W}}^{1}_{2}(\Omega)$.
\begin{equation}\label{eq2.28}
\|\boldsymbol{w}^{\varepsilon}\|_{2,\Omega}\leqslant\,
\varepsilon{M}_{\Omega}\|\mathbb{D}(x,\boldsymbol{w}^{\varepsilon})\|_{2,\Omega}.
\end{equation}
If $\Omega\subset\bigcup_{|\boldsymbol{k}|=1}^{n^{3}}\Omega^{\boldsymbol{k},\varepsilon}$
and $\boldsymbol{w}^{\varepsilon}\in\stackrel{\!\!\circ}{\mathbb{W}}^{1}_{2}
(\Omega^{\,\boldsymbol{k},\varepsilon})$, $\boldsymbol{k}=(k_{1},k_{2},k_{3}),\in\mathbb{Z}$, then
\begin{equation*}
\displaystyle\,\int_{\Omega^{\boldsymbol{k},\varepsilon}}|\boldsymbol{w}^{\varepsilon}|^{2}dx
\,\leqslant\,\varepsilon^{2}\,M_{\Omega}\int_{\Omega^{\boldsymbol{k},\varepsilon}}
|\mathbb{D}\big(x,\boldsymbol{w}^{\varepsilon}(\boldsymbol{x})\big)|^{2}dx
\end{equation*}
and
\begin{equation}\label{eq2.29}
\int_{\Omega}|\boldsymbol{w}^{\varepsilon}(\boldsymbol{x})|^{2}dx\,\leqslant\,
\varepsilon^{2}\,M_{\Omega}\int_{\Omega}
|\mathbb{D}\big(x,\boldsymbol{w}^{\varepsilon}(\boldsymbol{x})\big)|^{2}dx.
\end{equation}
\end{lemma}

\begin{remark}\label{Rm2.2}
A similar result
\begin{equation*}
\int_{\Omega}(|\boldsymbol{w}^{\varepsilon}(\boldsymbol{x})-
\boldsymbol{w}^{\varepsilon}_{\Omega}(\boldsymbol{x})|^{2}dx\,\leqslant\,
\varepsilon^{2}\,M_{\Omega}\int_{\Omega}
|\mathbb{D}\big(x,\boldsymbol{w}^{\varepsilon}(\boldsymbol{x})\big)|^{2}dx,
\end{equation*}
where $\boldsymbol{w}^{\varepsilon}_{\Omega}(\boldsymbol{x})=
\frac{1}{|\Omega|}\int_{\Omega}\boldsymbol{w}^{\varepsilon}(\boldsymbol{z})dz$ and $|\Omega|$
is a volume of the domain $\Omega$, called the \textbf{\emph{Poincar\'{e}-Wirtinger inequality}}, holds for any $\boldsymbol{w}^{\varepsilon}\in\mathbb{W}^{1}_{2}(\Omega)$. For the proof, see~\cite{32}.
\end{remark}

\subsection{The simplest embedding theorem}

\begin{lemma}\label{l2.2}
Let  $\Omega\subset\mathbb{R}^{3}$ with piecewise $\mathbb{C}^{1}$ boundary.

Then for any function ${w}\in\mathbb{W}^{1}_{2}(\Omega)$ ($\boldsymbol{w}\in\mathbb{W}^{1}_{2}(\Omega)$) identically equal zero on some part of the boundary $\partial\Omega$ with strictly positive surface measure holds true the estimate
\begin{equation}\label{eq2.30}
\|{w}\|_{2}\leqslant{M}\|\nabla{w}\|_{2,\Omega},\,\,
(\|{\boldsymbol{w}}\|_{2}\leqslant{M}\|\mathbb{D}(x,\boldsymbol{w})\|_{2,\Omega}).
\end{equation}
The constant ${M}$ is bounded for bounded $\Omega$.
\end{lemma} (See~\cite{33,34}).

\subsection{Extension lemma}

Extension results are very important in homogenization (Zhikov et al., see~\cite{10,11}).

For example, some sequence has different properties in different domains and only the properties of the sequence in the first domain permit one to choose a convergent subsequence. Therefore, we must preserve the best properties of the sequence and apply the extension from the first domain onto the second one.
Fortunately, all the indicated results apply to our case for a structure
with special periodicity because in each cell of periodicity $\Omega^{\,\boldsymbol{k},\varepsilon}$ we may directly use the method suggested in~\cite{10} for soft inclusions (see chapter 3 ``Elementary Soft and Stiff Problems'', section 3.1, pp. 86-95).

In the following lemma we consider solutions $\{\boldsymbol{w}^{\varepsilon}_{j},\,j=f,s\}$ to the problem $\mathbb{B}^{\varepsilon}(r)$.
\begin{lemma} \label{l2.3}
1) Let $\{\boldsymbol{w}^{\varepsilon}_{f}\}$ be a bounded sequence in $\mathbb{W}^{1,0}_{2}\big(\Omega^{\varepsilon}_{f,T}(r)\big)$.

Then for all $\varepsilon>0$ there exist extension operator
$\mathbb{E}_{f}:\mathbb{W}^{1,0}_{2}\big(\Omega_{f,T}(r)\big)\rightarrow\mathbb{W}^{1,0}_{2}(\Omega_{T})$ $\mathbb{E}_{f}({\boldsymbol{w}}_{f}^{\varepsilon})=
\widetilde{\boldsymbol{w}}^{\,\varepsilon}_{f}\in\mathbb{W}^{1,0}_{2}(\Omega_{T})$, such that $\widetilde{\boldsymbol{w}}_{f}^{\,\varepsilon}\in\mathbb{W}^{1,0}_{2}(\Omega_{T})$ 
\begin{multline*}
\|\widetilde{\boldsymbol{w}}_{f}^{\,\varepsilon}\|_{2,\Omega_{T}}\leqslant{M}
\|\boldsymbol{w}_{f}^{\,\varepsilon}\|_{2,\Omega_{f,T}(r)},\,\,\,
\|\mathbb{D}(x,\widetilde{\boldsymbol{w}}_{f}^{\,\varepsilon})\|_{2,\Omega_{T}}\leqslant{M}
\|\mathbb{D}(x,\boldsymbol{w}_{f}^{\,\varepsilon})\|_{2,\Omega_{f,T}(r)},
\\
\|\chi^{\,\varepsilon}(\varepsilon
\frac{\partial\widetilde{\boldsymbol{w}}_{f}^{\varepsilon}}{\partial{t}})\|_{2,\Omega_{T}}=
\|\chi^{\varepsilon}
(\varepsilon\frac{\partial\boldsymbol{w}_{f}^{\varepsilon}}{\partial{t}})\|_{2,\Omega_{f,T}}.
\end{multline*}
2) Let $\{\boldsymbol{w}^{\varepsilon}_{s}\}$ be a bounded sequence in $\mathbb{W}^{1,0}_{2}\big(\Omega^{\varepsilon}_{s,T}(r)\big)$. Due to the boundary condition $\widetilde{\boldsymbol{w}}_{f}^{\varepsilon}=0$ at $\Gamma^{\varepsilon}(r)$ we set $\boldsymbol{w}_{s}^{\,\varepsilon}=0$ at $\Gamma^{\varepsilon}(r)$ and $\boldsymbol{w}_{s}^{\,\varepsilon}\equiv0\,\in\Omega^{\varepsilon}_{f,T}(r)$, which immediately implies 
$\mathbb{D}(x,\widetilde{\boldsymbol{w}}_{s}^{\,\varepsilon})\equiv{0}$ in $\overline{\Omega}^{\,\varepsilon}_{f,T}(r)$. Coming back to the continuity condition (\ref{eq2.8}) at the free boundary  $\Gamma^{\varepsilon}(r)$ we obtain $\mathbb{D}(x,\widetilde{\boldsymbol{w}}_{f}^{\,\varepsilon})=0$ at the free boundary $\Gamma^{\varepsilon}(r)$.

Then for all $\varepsilon>0$ there exist extension operator
$\widetilde{\boldsymbol{w}}_{s}^{\,\varepsilon}=\mathbb{E}_{s}(\boldsymbol{w}^{\varepsilon}_{s})$, $\mathbb{E}_{s}: \mathbb{W}^{1,0}_{2,s}(\Omega_{T})\rightarrow\mathbb{W}^{1,0}_{2}(\Omega_{T})$, such that $\widetilde{\boldsymbol{w}}_{s}^{\,\varepsilon}=0$ in $\overline{\Omega}^{\,\varepsilon}_{f,T}$,
$\cap\stackrel{\!\!\circ}{\mathbb{W}}^{1,0}_{2}\big(\overline{\Omega}^{\,\varepsilon}_{f,T}(r)\big)$ and
\begin{multline}\label{eq2.31}
(\widetilde{\boldsymbol{w}}_{s}^{\,\varepsilon}-\boldsymbol{w}^{\,\varepsilon}_{s})
(1-\chi^{\,\varepsilon})=0,\,\,\big(\mathbb{D}(x,\widetilde{\boldsymbol{w}}_{s}^{\,\varepsilon})-
\mathbb{D}(x,\boldsymbol{w}_{s}^{\,\varepsilon})\big)(1-\chi^{\varepsilon})=0,
\\
\widetilde{\boldsymbol{w}}_{s}^{\,\varepsilon}(\boldsymbol{x},t)=0,\,\,
\boldsymbol{x}\in\overline{\Omega}^{\,\varepsilon}_{f}(r),\,\,
|\widetilde{\boldsymbol{w}}_{s}^{\,\varepsilon}|^{(1,0})_{\Omega_{T}}
\leqslant{M}\|(1-\chi^{\varepsilon})\boldsymbol{w}_{s}^{\,\varepsilon}\|^{(1,0})_{2,\Omega_{s,T}(r)},
\end{multline}
where $M$ is independent of $\varepsilon$.

Moreover,
\begin{equation}\label{eq2.32}
\max_{0<t<T}|\widetilde{\boldsymbol{w}}_{s}^{\,\varepsilon}(.,t)|_{\Omega}^{(2+\alpha)}\leqslant\,M.
\end{equation}
\end{lemma}

\begin{proof}
To  prove the first statement we note that there are several options for the extensions of $\boldsymbol{w}^{\,\varepsilon}_{f}$. We choose the extension
\begin{equation}\label{eq2.33}
\widetilde{\boldsymbol{w}}^{\,\varepsilon}_{f}=
\chi^{\,\varepsilon}\boldsymbol{w}^{\,\varepsilon}_{f}-
(1-\chi^{\,\varepsilon})\boldsymbol{w}_{s}^{\varepsilon}
\end{equation}
for which
\begin{equation*}
\widetilde{\boldsymbol{w}}^{\,\varepsilon}_{f}(\boldsymbol{x},t)=0,\,\,
\boldsymbol{x}\in\overline{\Omega}^{\,\varepsilon}_{s,T}(r),\,\,
\chi^{\varepsilon}\mathbb{D}(x,\widetilde{\boldsymbol{w}}_{f}^{\,\varepsilon})=
\mathbb{D}(x,\boldsymbol{w}_{f}^{\,\varepsilon})\,\,\,(\boldsymbol{x},t)\in\overline{\Omega}^{\,\varepsilon}_{s,T}(r).
\end{equation*}
Thus, to prove the statement we just put $\widetilde{\boldsymbol{w}}^{\,\varepsilon}_{f}(\boldsymbol{x},t)=0$ in $\overline{\Omega}_{s,T}$.

It is easy to see that $\widetilde{\boldsymbol{w}}_{s}^{\,\varepsilon}$ satisfies all conditions of the lemma.

The estimate (\ref{eq2.32}) follows from Theorem 1.1, Chapter 3, \S 1 in~\cite{29}.
\end{proof}
The following lemma is true:
\begin{lemma} \label{l2.4}
Let $\{{c}^{\,\varepsilon}\}$ be bounded sequence in $\mathbb{W}^{1,0}_{2}\big(\Omega^{\varepsilon}_{f,T}(r)\big)$ and ${c}_{0}\in\mathbb{W}^{1,0}_{2}\big(\Omega^{\varepsilon}_{f,T}(r)\big)$.

1) Then for all $\varepsilon>0$ there exist extensions $\widetilde{c}^{\,\varepsilon}$, such that
\begin{equation}\label{eq2.34}
\|(\widetilde{c}^{\,\varepsilon}-{c}_{0})\|_{2,\Omega_{T}}+
\|\nabla(\widetilde{c}^{\,\varepsilon}-{c}_{0})\|_{2,\Omega_{T}}\leqslant{M},
\end{equation}
\end{lemma}
For proof see~\cite{26,27}.

Let $\widetilde{\mathbb{P}}^{\,\varepsilon}=
\chi^{\varepsilon}\mathbb{P}_{f}+(1-\chi^{\varepsilon})\mathbb{P}_{f}$.

\begin{remark}\label{2.3}
Due to the choice of function $p^{0}(\boldsymbol{x})$
\begin{equation*}
\widetilde{\mathbb{P}}^{\,\varepsilon}<\boldsymbol{n}>=0
\end{equation*}
at the boundary ${S}^{1}\cup{S}^{2}$. Here $\boldsymbol{n}$ is the normal vector to ${S}^{1}$ and ${S}^{2}$.
\end{remark}

\subsection{H\"older's inequality}

\begin{lemma} \label{l2.5}
For any $f,\,g\in\mathbb{L}_{2}(\Omega)$ holds true the H\"older's inequality
\begin{equation}\label{eq2.35}
\|f\,g\|_{1,\Omega}\leqslant\|f\|_{2,\Omega}\|g\|_{2,\Omega}.
\end{equation}
\end{lemma}
For details, see~\cite{34}.

\subsection{Two-scale convergence method}

In the present section we consider 1-periodic in the variable $\boldsymbol{y}\in{\textbf{Y}}$ functions  ${W}(\boldsymbol{y};\boldsymbol{x},t)$ with $(\boldsymbol{x},t)\in \Omega_{T}$.

\begin{definition}\label{de2.1}
\emph{The sequence $\{w^{\,\varepsilon}\}\subset \mathbb{L}_{2}(\Omega_{T})$,
is said to be two -- scale convergent to the function
$W(\boldsymbol{y};\boldsymbol{x},t)\in \mathbb{L}_{2}(\Omega_{T}\times Y)$,
which is 1 -- periodic in the variable $\boldsymbol{y}\in Y$ ( notation
$w^{\,\varepsilon}\,\stackrel{\!\!2-sc.}{\rightarrow}\,W(\boldsymbol{x},t;\boldsymbol{y})$),
if for any smooth function $\sigma=\sigma(\boldsymbol{y};\boldsymbol{x},t)$,
1 -- periodic in the variable $\boldsymbol{y}$ is valid the equality}
\begin{multline}\label{eq2.36}
\,\,\,\,\,\,\,\,\,\,\,\,\,\,\,\,\,\,\,\,
\lim_{{\varepsilon}\to 0}\int\int_{\Omega_{T}} w^{\,\varepsilon}(\boldsymbol{x},t) \sigma(\frac{\boldsymbol{x}}{\varepsilon};\boldsymbol{x},t)dxdt=
\\
\int\int_{\Omega_{T}}\big(\int_Y{W}(\boldsymbol{y};\boldsymbol{x},t)
\sigma(\boldsymbol{x},t;\boldsymbol{y})dy\big)dxdt.
\,\,\,\,\,\,\,\,\,\,\,\,\,\,\,\,\,\,\,\,\,\,\,\,\,\,\,\,\,\,\,\,\,\,\,\,\,\,\,\,
\end{multline}
\end{definition}
Note that weak and  two--scale convergence are connected by the relation:
\begin{multline*}
\mbox{if}\,\,{u}^{\,\varepsilon}\,\stackrel{\!\! 2-sc.}{\rightarrow}\,
{U}(\boldsymbol{y};\boldsymbol{x},t)\,\,\,\,\mbox{(converges two -- scale)},
\\
\mbox{then}\,\,{u}^{\,\varepsilon}(\boldsymbol{x},t)
\rightharpoonup \int_Y {U}(\boldsymbol{y};\boldsymbol{x},t) d\boldsymbol{y}=u(\boldsymbol{x},t)
\,\,\mbox{(converges weakly)}.
\end{multline*}
The existence and basic properties of two--scale convergent sequences are proved in the following theorem:

\begin{theorem} \label{t2.1} (Nguetseng's Theorem~\cite{14})

\textbf{1.} Any bounded in $\mathbb{L}_{2}(\Omega_{T})$ sequence $\{\boldsymbol{w}^{\,\varepsilon}\}$
contains some subsequence two -- scale convergent to some function
\\
$\boldsymbol{W}(\boldsymbol{y};\boldsymbol{x},t)$, $\boldsymbol{W}\in\mathbb{L}_{2}
(\Omega_{T}\times \textbf{Y})$, 1 -- periodic in the variable $\boldsymbol{y}$.

\textbf{2.} Let sequences $\{\boldsymbol{w}^{\,\varepsilon}\}$ and $\{\varepsilon\mathbb{D}(x,\boldsymbol{w}^{\,\varepsilon})\}$ be uniformly bounded in $\mathbb{L}_{2}(\Omega_{T})$.

Then there exists the function $\boldsymbol{W}=\boldsymbol{W}(\boldsymbol{y};\boldsymbol{x},t)$, 1 -- periodic
in $\boldsymbol{y}$, and the sequence $\{\boldsymbol{w}^{\,\varepsilon}\}$ such that
$\boldsymbol{W},\,\mathbb{D}(y,\boldsymbol{W})\in \mathbb{L}_{2}(\Omega_{T}\times\textbf{Y})$,
and sequences $\{\boldsymbol{w}^{\,\varepsilon}\}$ and $\{\varepsilon\mathbb{D}(x,\boldsymbol{w}^{\,\varepsilon})\}$ (for simplicity we keep the same indices for subsequences) two -- scale converge in
$\mathbb{L}_{2}(\Omega_{T})$ to $\boldsymbol{W}$ and $\mathbb{D}(y,\boldsymbol{W})$ correspondingly.

\textbf{3.} Let sequences $\{\boldsymbol{w}^{\,\varepsilon}\}$ and
$\{D(x,\boldsymbol{w}^{\,\varepsilon})\}$ be bounded in $\mathbb{L}_{2}(\Omega_{T})$.

Then there are functions $\boldsymbol{w}(\boldsymbol{x},t),\,\boldsymbol{w}\in \mathbb{W}^{1,0}_{2}(\Omega_{T})$ and
$\boldsymbol{W}(\boldsymbol{y};\boldsymbol{x},t),\,\boldsymbol{W}\in \mathbb{L}_{2}(\Omega_{T}\times \textbf{Y})\cap\mathbb{W}^{1,0}_{2}(\textbf{Y})$,
subsequence from $\{\mathbb{D}(x,\boldsymbol{w}^{\,\varepsilon})\}$ such that the function
$\boldsymbol{W}$ is 1 -- periodic in $\boldsymbol{y}$, $\mathbb{D}(x,\boldsymbol{w}) \in \mathbb{L}_{2}(\Omega_{T})$,
$D(y,\boldsymbol{W}) \in \mathbb{L}_{2}(\Omega_{T}\times\textbf{Y})$, and the sequence
$\{\mathbb{D}(x,\boldsymbol{w}^{\,\varepsilon})\}$ two -- scale converges to the function
$\mathbb{D}(x,\boldsymbol{w})+D(y,\boldsymbol{W})$.
\end{theorem}

\subsection{Strong convergence criteria in $\mathbb{L}_{2}(\Omega)$}

\begin{definition}\label{de2.2}
We say that the function $\widetilde{c}(\boldsymbol{x},t)$, bounded in $\mathbb{L}_{2}(\Omega_{T})$, possesses a time derivative
$\displaystyle\,\frac{\partial\widetilde{c}}{\partial{t}}\in\mathbb{L}_{2}\big(0,T;\mathbb{W}^{-1}_{2}(\Omega)\big)$, if
\begin{equation*}
|\int\int_{\Omega_{T}}\widetilde{c}\frac{\partial \xi}{\partial t}dxdt|\leqslant \,M_{u}\,|\int\int_{\Omega_{T}} |\nabla\xi|^{2}dxdt|^{\frac{1}{2}}
\end{equation*}
for all functions $\xi\in\mathbb{W}^{1,1}_{2}(\Omega_{T})$ with some positive constant $M_{u}$ independent of $\xi$.
\end{definition}
\begin{lemma} \label{l2.6}
Let sequences $\{\widetilde{c}^{\,\varepsilon}\}$ and $\{\nabla\widetilde{c}^{\,\varepsilon}\}$ be uniformly bounded in the space $\mathbb{L}_{2}(\Omega_{T})$, and the sequence of derivatives
$\{\frac{\partial\widetilde{c}^{\,\varepsilon}}{\partial{t}}\}$
are uniformly bounded in the space $\mathbb{L}_{2}\big(0,T;\mathbb{W}^{-1}_{2}(\Omega)\big)$.

Then there exists a subsequence of the sequence $\{\widetilde{c}^{\,\varepsilon}\}$ that strongly converges in $\mathbb{L}_{2}(\Omega_{T})$.
\end{lemma}
For the proof, see~\cite{30}.

The generalization of this lemma for a periodic structure with characteristic function
$\chi^{\,\varepsilon}(\boldsymbol{x})=\chi(\boldsymbol{x};\frac{\boldsymbol{x}}{\varepsilon})$
has been proved by A.~Meirmanov and O.~Galtsev in~\cite{35}.

\begin{lemma} \label{l2.7}
Let $\chi^{\varepsilon}(\boldsymbol{x})=\chi(\boldsymbol{x},\frac{\boldsymbol{x}}{\varepsilon})$,
where $\chi(\boldsymbol{x},\boldsymbol{y})$ is a 1--periodic in $\boldsymbol{y}$ function, the sequences
$\{\widetilde{c}^{\,\varepsilon}\}$ and $\{\nabla\widetilde{c}^{\,\varepsilon}\}$ are uniformly bounded in $\mathbb{L}_{2}(\Omega_{T})$,
and the sequence $\{\chi^{\varepsilon}\frac{\partial\widetilde{c}^{\,\varepsilon}}{\partial{t}}\}$
is uniformly bounded in $\mathbb{L}_{2}\big(0,T;\mathbb{W}^{-1}_{2}(\Omega)\big)$.

Then there exists some subsequence of $\{\widetilde{c}^{\,\varepsilon}\}$ that converges strongly in $\mathbb{L}_{2}(\Omega_{T})$.
\end{lemma}
\begin{remark}\label{Rm2.4}
We denote the norm of an element $\varphi\in\mathbb{L}_{2}\big(0,T;\mathbb{W}^{-1}_{2}(\Omega)\big)$ as $\|\varphi\|_{W^{-1}_{2}}$.
\end{remark}
We further formulate a similar result for  structures with special periodicity for pore and solid skeleton structures.

\begin{theorem}\label{t2.2}

Let the structure $\chi(r;\boldsymbol{y})$ of the pore space be given by formula \eqref{eq2.25}, where
$r\in \mathfrak{M}_{T}$ and  
\begin{multline}\label{eq2.37}
\mathfrak{M}_{(0,T)}=\{r\,\in \mathbb{H}^{2+\gamma,\frac{2+\gamma}{2}}(\overline{\Omega}_{T}),\,\,
0\,\leqslant\,r(\boldsymbol{x},t)\leqslant\,\frac{1}{2},\,\,\,
-\theta\,\leqslant\,\frac{\partial r}{\partial t}(\boldsymbol{x},t)\leqslant\,0,\\
0<\gamma<1,\,\,\,\theta=\mbox{const}>0;\,\,\,
|r|^{(2+\gamma)}_{\Omega_{T}}\leqslant\,M_{0}\}.
\end{multline}
Then any sequence $\{\widetilde{c}^{\,\varepsilon}\}$, such that
\begin{equation*}
\|\widetilde{c}^{\,\varepsilon}\|_{2,\Omega_{T}}+
\|\nabla\widetilde{c}^{\,\varepsilon}\|_{2,\Omega_{T}}+
\|\frac{\partial}{\partial t}\widetilde{c}^{\,\varepsilon}\|_{W^{-1}_{2}}\,\leqslant\,M,
\end{equation*}
where $M$ does not depend on $\varepsilon$, contains a subsequence, strongly convergent in $\mathbb{L}_{2}(\Omega_{T})$.
\end{theorem}
For the proof, see Theorem~2.2 in~\cite{18}.

\subsection{Equivalent formulation of the boundary condition \eqref{eq2.11} as an integral identity}

\begin{lemma}\label{l2.8}
Under the condition
\begin{equation*}
\alpha^{\varepsilon}=\varepsilon\,\theta,\,\,\,\beta^{\varepsilon}=\varepsilon,
\end{equation*}
where $\theta$ is a given positive constant, the boundary condition \eqref{eq2.11} is equivalent to the integral identity
\begin{equation}\label{eq2.38}
\int_{0}^{t_{0}}\int_{\Omega_{}}\chi^{\varepsilon}(\boldsymbol{x},t)\Big(-\frac{\partial}{\partial t}
\big((\zeta\,\boldsymbol{a}^{\varepsilon})\cdot\boldsymbol{\xi}_{0}^{\varepsilon})\big)+
\varepsilon\,\nabla\cdot(\zeta(\widetilde{c}^{\,\varepsilon}-{c}^{0})\boldsymbol{\xi}_{0}^{\varepsilon})\Big)dxdt=0
\end{equation}
which is valid for any smooth functions
$\boldsymbol{\xi}_{c}^{\,\varepsilon}(r,\boldsymbol{x})=
\boldsymbol{\xi}_{c}(r,\frac{\boldsymbol{x}}{\delta})$, functions $\zeta$, vanishing at $t=0$ and at $t=t_{0}$ and at boundary $\partial\Omega$, and functions
$\boldsymbol{a}_{c}^{\,\varepsilon}(r,\boldsymbol{x})=
\boldsymbol{a}_{c}(r,\frac{\boldsymbol{x}}{\delta})$,
such that $\boldsymbol{a}_{c}$ vanishes outside of some small neighbourhood of $\gamma_{c}(r)$ and $\boldsymbol{a}_{c}(r,\boldsymbol{y})=\boldsymbol{n}_{c}(r)$, where $\boldsymbol{n}_{c}(r)$
is the unit normal to the surface $\gamma_{c}(r)=\{\boldsymbol{y}\in\textbf{Y}: |\boldsymbol{y}|=r\}$,
outward to the domain $\textbf{Y}_{f}(r)$.
\end{lemma}
For the proof, see Lemma 4.2 in~\cite{18}.

\subsection{Equivalent formulation of the dynamic problem $\mathbb{A}^{\varepsilon}$ 
as an integral identity}

We assume that $p^{0}(\boldsymbol{x})$ is given bounded function, $p^{0}\in\mathbb{C}^{1}(\overline{\Omega})$ and $p^{0}(\boldsymbol{x})=p^{j}$=\emph{const} for $\boldsymbol{x}\in{S}^{j}$,\,j=1,2.
\begin{definition}\label{d3}
Let the structure $\chi^{*,\varepsilon}$ of the pore space $\Omega^{\varepsilon}_{f,T}(r^{*})$ be given by the function $r^{*}\in \mathfrak{M}_{(0,T)}$, ${p}^{\varepsilon}=
\chi^{*,\varepsilon}({p}^{\varepsilon}_{f}-{p}^{0})+(1-\chi^{*,\varepsilon})({p}^{\varepsilon}_{s}-{p}^{0})$,
$\mathbb{P}^{\varepsilon}=\chi^{*,\varepsilon}\mathbb{P}^{\varepsilon}_{f}+
(1-\chi^{*,\varepsilon})\mathbb{P}^{\varepsilon}_{s}$,
$\mathbb{P}^{\varepsilon}_{f}=\varepsilon^{2}\mu_{1}
\mathbb{D}(x,\varepsilon\frac{\partial\boldsymbol{w}_{f}^{\varepsilon}}{\partial{t}})-
({p}^{\varepsilon}_{f}-{p}^{0})\mathbb{I}$, $\mathbb{P}^{\varepsilon}_{s}=
\lambda_{0}\mathbb{D}(x,\boldsymbol{w}_{s}^{\varepsilon})-({p}^{\varepsilon}_{s}-{p}^{0})\mathbb{I}$.

We say that functions
$\boldsymbol{w}^{\varepsilon}_{f}\in\mathbb{W}^{1,0}_{2}\big(\Omega_{f,T}(r^{*})\big)$,
$\boldsymbol{w}^{\varepsilon}_{s}\in\mathbb{W}^{1,0}_{2}\big(\Omega_{s,T}(r^{*})\big)$,
${p}^{\varepsilon}_{f}$ and ${p}^{\varepsilon}_{f}$ are a weak solution to the Dynamic problem $\mathbb{A}^{\varepsilon}$
if the following equations hold: (\ref{eq2.2}), (\ref{eq2.4}) and integral identity
\begin{multline}\label{eq2.39}
-\int_{0}^{t_{0}}\int_{\Omega}(\nabla{p}^{0}\cdot\boldsymbol{\varphi})dxdt=
\\
\int_{0}^{t_{0}}\int_{\Omega}\Big(\chi^{*,\varepsilon}\varepsilon^{2}\mu_{1}
\mathbb{D}(x,\frac{\partial\boldsymbol{w}_{f}^{\varepsilon}}{\partial{t}})+
(1-\chi^{*,\varepsilon})\lambda_{0}\mathbb{D}(x,\boldsymbol{w}_{s}^{\varepsilon})-
\\
\big(\chi^{*,\varepsilon}({p}^{\varepsilon}_{f}-{p}^{0})+(1-\chi^{*,\varepsilon})
({p}^{\varepsilon}_{s}-{p}^{0})\mathbb{I}\big)\Big):\mathbb{D}(x,\boldsymbol{\varphi})dxdt
\end{multline}
for any arbitrary smooth functions $\boldsymbol{\varphi}$,  vanishing at the boundary
$\big(S^{1}\cup S^{2}\big)\times(0,T)$ and satisfying conditions
\begin{multline}\label{eq2.40}
\boldsymbol{\varphi}(\boldsymbol{x}_{0}+0)=\boldsymbol{\varphi}(\boldsymbol{x}_{0}-0),\,\,\,
\boldsymbol{x}_{0}\in\Gamma^{\varepsilon}(r^{*})
\\
\,\,\,\,\,\,\,\,\,\boldsymbol{\varphi}(\boldsymbol{x}_{0}+0)=
\lim_{\boldsymbol{x}\rightarrow\boldsymbol{x}_{0}}\boldsymbol{\varphi}(\boldsymbol{x}),
\,\,\boldsymbol{x}\in\Omega^{\varepsilon}_{f,t_{0}}(r^{*}),\,\,\,
\boldsymbol{x}_{0}\in\Gamma^{\varepsilon}(r^{*}),
\\
\boldsymbol{\varphi}(\boldsymbol{x}_{0}-0)=
\lim_{\boldsymbol{x}\rightarrow\boldsymbol{x}_{0}}
\boldsymbol{\varphi}(\boldsymbol{x}),\,\,\boldsymbol{x}\in \Omega^{\varepsilon}_{s,t_{0}}(r^{*}),\,\,\,
\boldsymbol{x}_{0}\in \Gamma^{\varepsilon}(r^{*})
\end{multline}
at the free boundary $\Gamma^{\varepsilon}(r^{*})$.
\end{definition}

\subsection{Equivalent formulation of the diffusion problem $\mathbb{A}^{\varepsilon}$ 
as an integral identity}

\begin{definition}\label{d4}
Let the structure $\chi^{*,\varepsilon}$ of the pore space $\Omega^{\,\varepsilon}_{f,T}(r^{*})$ be given by the function $r^{*}\in \mathfrak{M}_{(0,T)}$.

We say that function ${c}^{\varepsilon}$ is a weak solution to the diffusion problem $\mathbb{A}^{\varepsilon}$,
if the integral identity holds
\begin{multline}\label{eq2.41}
\int_{\Omega}\chi^{*,\varepsilon}(.,{t}_{0})
\big({c}^{\varepsilon}(.,{t}_{0})+
\frac{\beta^{\varepsilon}}{\alpha^{\varepsilon}}\big)\xi\chi^{*,\varepsilon}(.,{t}_{0})dx-
\int_{\Omega}\chi^{*,\varepsilon}(.,0)\big({c}^{0}+
\frac{\beta^{\varepsilon}}{\alpha^{\varepsilon}}\big)\,\xi\chi^{*,\varepsilon}(.,0)\,dx+
\\
\int_{0}^{t_{0}}\int_{\Omega}\chi^{*,\varepsilon}\big(-({c}^{\varepsilon}+
\frac{\beta^{\varepsilon}}{\alpha^{\varepsilon}})\,
\frac{\partial\xi}{\partial{t}}+\nabla\xi\cdot(\alpha_{c}
\nabla{c}^{\varepsilon}\big)dxdt=0
\end{multline}
for any arbitrary smooth functions $\xi$, vanishing at the boundary $\big(S^{1}\cup S^{2}\big)\times(0,T)$.
\end{definition}

\begin{remark}\label{Rm5}
In deriving the integral identity (\ref{eq2.41}), we used the boundary condition (\ref{eq2.11}) on the free boundary, so that the term containing the integral over this boundary vanishes.
\end{remark}

\subsection{Equivalent formulation of the dynamic problem $\mathbb{B}^{\varepsilon}(r)$ 
as integral identities}

By $\widetilde{\pi}^{\varepsilon}$ we denote the \textbf{\emph{Antiderivative}} of the function $\widetilde{p}^{\varepsilon}_{f}$:
$\widetilde{p}^{\,\varepsilon}_{f}(\boldsymbol{x},t)=
\int_{0}^{t}\widetilde{\pi}^{\varepsilon}(\boldsymbol{x},\tau)d\tau$
\begin{definition}\label{df5}
Let the structure $\chi^{\varepsilon}$ of the pore space $\Omega^{\varepsilon}_{f,T}(r)$ be given by the function $r\in \mathfrak{M}_{(0,T)}$ and
\begin{multline*}
I_{1}^{\varepsilon}(\boldsymbol{\varphi})=\int_{0}^{t_{0}}\int_{\Omega}\chi^{\varepsilon}
\big(\mu_{1}\varepsilon^{2}\mathbb{D}(x,\frac{\partial\boldsymbol{w}^{\varepsilon}_{f}}{\partial t})\big):
\mathbb{D}(x,\boldsymbol{\varphi})dxdt,
\\
I_{2}^{\varepsilon}(\boldsymbol{\varphi})=\int_{0}^{t_{0}}\int_{\Omega}\lambda_{0}(1-\chi^{\varepsilon})
\big(\mathbb{D}(x,\boldsymbol{w}_{s}^{\varepsilon})+
{c}_{s}^{2}\nabla\cdot\boldsymbol{w}^{\varepsilon}_{s}\big):\mathbb{D}(x,\boldsymbol{\varphi})dxdt,
\\
I_{3}^{\varepsilon}(\boldsymbol{\varphi})=\varepsilon\int_{0}^{t_{0}}\int_{\Omega}\chi^{\varepsilon}
(\frac{\partial\boldsymbol{w}^{\varepsilon}_{f}}{\partial t}\cdot\boldsymbol{\varphi})dxdt.
\end{multline*}
We say that functions
$\widetilde{\boldsymbol{w}}^{\varepsilon}_{f},\,\widetilde{\boldsymbol{w}}^{\varepsilon}_{s}
\in\mathbb{W}^{1,0}_{2}(\Omega_{T})$,
$\sqrt{\varepsilon}\frac{\partial\widetilde{\boldsymbol{w}}_{f}^{\varepsilon}}{\partial{t}},\,
\widetilde{{p}}^{\varepsilon},\,\frac{\partial{\widetilde{\pi}}^{\varepsilon}}{\partial{t}}
\in\mathbb{L}_{2}\big(\Omega_{T}$ are a weak solution to the Dynamic problem $\mathbb{B}^{\varepsilon}(r)$, if they satisfy continuity equation (\ref{eq2.2}), boundary condition (\ref{eq2.14}) and integral identity
\begin{equation}\label{eq2.42}
-\int_{0}^{t_{0}}\int_{\Omega}\chi^{\varepsilon}(\nabla{p}^{0}+
\varepsilon\frac{\partial\widetilde{\boldsymbol{w}}_{f}^{\varepsilon}}{\partial{t}})
\cdot\boldsymbol{\varphi}dxdt=
I_{1}^{\varepsilon}(\boldsymbol{\varphi})+I_{2}^{\varepsilon}(\boldsymbol{\varphi})+
I_{3}^{\varepsilon}(\boldsymbol{\varphi})
\end{equation}
with test function $\boldsymbol{\varphi}$, vanishing at the boundary ${S}^{0}\times(0,T)$
and satisfying condition (\ref{eq2.40}) at the boundary $\Gamma^{\varepsilon}(r)$.
\end{definition}

\begin{definition}\label{Df7}
Let the structure $\chi^{\varepsilon}$ of the pore space $\Omega^{\varepsilon}_{f,T}(r)$ be given by the function $r\in \mathfrak{M}_{(0,T)}$.

We say that functions
$\widetilde{\boldsymbol{w}}^{\,\varepsilon}_{s}\in\mathbb{W}^{1,0}_{2}
\big(\Omega_{f,T}(r)\big)$, $\widetilde{\,{p}}^{\varepsilon}_{s}
\in\mathbb{L}_{2}\big(\Omega_{T}(r)\big)$ are a weak solution to the Dynamic problem $\mathbb{B}^{\varepsilon}(r)$ for the solid component, if continuity equation (\ref{eq2.4}) and integral identity
\begin{equation}\label{eq2.43}
\int_{0}^{t_{0}}\int_{\Omega}\Big(\big((1-\chi^{\varepsilon})(\nabla{p}^{0}\cdot\boldsymbol{\varphi})+
\lambda_{0}\mathbb{D}(x,\widetilde{\boldsymbol{w}}_{s}^{\varepsilon})+
{c}_{s}^{2}(\nabla\cdot\widetilde{\boldsymbol{w}}_{s}^{\varepsilon})
\mathbb{I}\big):\mathbb{D}(x,\boldsymbol{\varphi})\Big)dxdt=0
\end{equation}
for the solid component with arbitrary test function $\boldsymbol{\varphi}$, satisfying conditions (\ref{eq2.40}) at the boundary $\Gamma^{\varepsilon}(r)$, hold.
\end{definition}

In the identity (\ref{eq2.43}) we used the continuity equation (\ref{eq2.4}).

\begin{definition}\label{Df7b}
Let the structure $\chi^{\varepsilon}_{s}$ of the domain $\Omega^{\varepsilon}_{s,T}$ the pore space $\Omega^{\varepsilon}_{f,T}$ be given by the function $r\in \mathfrak{M}_{(0,T)}$.

We say that functions
$\widetilde{\boldsymbol{w}}^{\,\varepsilon}_{f}\in\mathbb{W}^{1,0}_{2}
\big(\Omega_{f,T}(r)\big)$, $\widetilde{\,{p}}^{\varepsilon}_{f}
\in\mathbb{L}_{2}\big(\Omega_{T}(r)\big)$ are a weak solution to the Dynamic problem $\mathbb{B}^{\varepsilon}(r)$ for the liquid component, if continuity equation (\ref{eq2.2})
and integral identity
\begin{multline}\label{eq2.44}
0=\int_{0}^{t_{0}}\int_{\Omega}\chi^{\varepsilon}\Big((\nabla{p}^{0}+
\varepsilon\frac{\partial\widetilde{\boldsymbol{w}}_{f}^{\varepsilon}}{\partial{t}})
\cdot\boldsymbol{\varphi}\big)+\varepsilon^{2}\mu_{1}
\mathbb{D}(x,\frac{\partial\widetilde{\boldsymbol{w}}_{f}^{\varepsilon}}{\partial{t}})\Big):
\mathbb{D}(x,\boldsymbol{\varphi})dxdt=
\\
-\int_{0}^{t_{0}}\int_{\Omega}\chi^{\varepsilon}\Big(\big((\nabla{p}^{0}t+
\varepsilon\widetilde{\boldsymbol{w}}_{f}^{\varepsilon})\cdot
\frac{\partial\boldsymbol{\varphi}}{\partial{t}}\big)dxdt-
\\
\int_{0}^{t_{0}}\int_{\Omega}\chi^{\varepsilon}\big(\varepsilon^{2}\mu_{1}
\mathbb{D}(x,\widetilde{\boldsymbol{w}}_{f}^{\varepsilon})\big):
\mathbb{D}(x,\frac{\partial\boldsymbol{\varphi}}{\partial{t}})\Big)dxdt
\end{multline}
for the liquid component with arbitrary test function $\boldsymbol{\varphi}$, satisfying conditions (\ref{eq2.40}) at the boundary $\Gamma^{\varepsilon}(r)$, hold.
\end{definition}

\subsection{Equivalent formulation of the diffusion problem $\mathbb{B}^{\varepsilon}(r)$ 
as an integral identity}  

\begin{definition}\label{d8}

Let the structure $\chi^{\varepsilon}$ of the pore space $\Omega^{\,\varepsilon}_{f,T}(r)$ be given by the function $r\in\mathfrak{M}_{(0,T)}$.

We say that function ${\widetilde{c}}^{\,\varepsilon}\in\mathbb{W}^{1,0}_{2}
\big(\Omega_{f,T}(r)\big)$ is a weak solution to the diffusion problem $\mathbb{B}^{\varepsilon}(r)$, if the integral identity holds
\begin{multline}\label{eq2.45}
\int_{\Omega}\chi^{\varepsilon}(.,{t}_{0})
\big({\widetilde{c}}^{\,\varepsilon}(.,{t}_{0})+
\frac{\beta^{\varepsilon}}{\alpha^{\varepsilon}}\big)\xi\chi^{\varepsilon}(.,{t}_{0})dx-
\\
\int_{0}^{t_{0}}\int_{\Omega}\chi^{\varepsilon}\big(-({\widetilde{c}}^{\,\varepsilon}-{c}^{0}+
\frac{\beta^{\varepsilon}}{\alpha^{\varepsilon}})\,\frac{\partial\xi}{\partial{t}}+
\nabla\xi\cdot(\alpha_{c}\nabla({\widetilde{c}}^{\,\varepsilon}-{c}^{0})\big)dxdt=0
\end{multline}
for any arbitrary smooth functions $\xi$, vanishing at the boundary $\big(S^{1}\cup S^{2}\big)\times(0,T)$
and at $t=0$.
\end{definition}

\subsection{Formal homogenization of the problem $\mathbb{A}^{\varepsilon}$}
As
\begin{equation*}
{\pi}(\boldsymbol{x},t)=
\lim_{\varepsilon\rightarrow0}{\widetilde{\pi}}^{\,\varepsilon}(\boldsymbol{x},t)=
\lim_{\varepsilon\rightarrow0}\int_{0}^{t}\widetilde{p}^{\,\varepsilon}_{f}(\boldsymbol{x},\tau)d\tau=
\int_{0}^{t}{p}_{f}(\boldsymbol{x},\tau)d\tau
\end{equation*}
we denote an antiderivative of the pressure ${p_{f}}$.
\begin{lemma} \label{l2.9}
Under the conditions of Theorem~\ref{t3.1} the formal homogenization $\mathbb{H}$ of the problem $\mathbb{A}^{\varepsilon}$ consisting of Darcy law of filtration
\begin{equation}\label{eq2.46}
\boldsymbol{w}_{f}=-\frac{1}{\mu_{1}}(B)^{(w)}(r)<\nabla{(\pi}-{p}^{0}t)>,\,\,
\nabla_{x}\cdot\boldsymbol{w}_{f}=0
\end{equation}
for the liquid displacements $\boldsymbol{w}_{f}$ and  the antiderivative ${\pi}$  of the liquid pressure ${p}_{f}$ in the domain $\Omega_{T}$, homogenized Lam\'{e} system
\begin{equation}\label{eq2.47}
\nabla\cdot\big(\lambda_{0}\mathfrak{N}^{(s)}:\mathbb{D}(x,\boldsymbol{w}_{s})+
{c}_{s}^{2}(\nabla\cdot\boldsymbol{w}_{s})\mathbb{I}\big)=\nabla{p}_{0},
\end{equation}
\begin{equation}\label{eq2.48}
\frac{1}{c^{2}_{s}}({p}^{\varepsilon}_{s}-{p}^{0})+\boldsymbol{w}_{s}=0
\end{equation}
for the solid displacements $\boldsymbol{w}_{s}$ and solid pressure $p_{s}$ and homogenized system, describing diffusion of the acid
\begin{equation}\label{eq2.49}
\frac{\partial}{\partial{t}}\big(m(r){c}\big)=
\nabla\cdot(\alpha_{c}(B)^{(c)}(r)<\nabla({c}-{c}^{0})>
\end{equation}
in the domain $\Omega_{T}$.

Differential equations are completed with the boundary and initial conditions
\begin{equation}\label{eq2.50}
{\pi}(\boldsymbol{x},t)={p}^{0}(\boldsymbol{x})t,\,\,\,\boldsymbol{x}\in{S}^{1}\cup{S}^{2},\,\,\,0<t<T,
\end{equation}
\begin{equation}\label{eq2.55}
\boldsymbol{w}_{f}\cdot\boldsymbol{n}=0,\,\,\,\boldsymbol{x}\in{S}^{0},\,\,\,0<t<T,
\end{equation}
\begin{equation}\label{eq2.51}
\boldsymbol{w}_{s}(\boldsymbol{x},t)=0,\,\,\,\boldsymbol{x}\in\partial\Omega,\,\,\,0<t<T
\end{equation}
\begin{equation}\label{eq2.52}
\big(\lambda_{0}\mathfrak{N}^{(s)}:\mathbb{D}(x,\boldsymbol{w}_{s})+
{c}_{s}^{2}\nabla\cdot\boldsymbol{w}_{s})\mathbb{I}\big)<\boldsymbol{n}>=0,\,\,
\boldsymbol{x}\in{S}^{1}\cup{S}^{2},\,\,0<t<T,
\end{equation}
\begin{equation}\label{eq2.53}
{c}(\boldsymbol{x},{t})={c}^{0}(\boldsymbol{x}),\,\,\,\boldsymbol{x}\in{S}^{1}\cup{S}^{2},\,\,t>0,
\end{equation}
\begin{equation}\label{eq2.54}
\frac{\partial{c}}{\partial{n}}(\boldsymbol{x},t)=0,\,\,\boldsymbol{x}\in S^{0},\,\,t>0,
\end{equation}
\begin{equation}\label{eq2.56}
{c}(\boldsymbol{x},0)={c}^{0}(\boldsymbol{x}),\,\,\,\boldsymbol{x}\in\Omega.
\end{equation}
In (\ref{eq2.46})-(\ref{eq2.56}), $\boldsymbol{n}$ is a normal vector to the boundary ${S}^{1}\cup{S}^{2}$, the symmetric strictly positive definite constant matrices $(B)^{({w})}(r)$ and $(B)^{({c})}(r)$ and the symmetric strictly positive definite fourth-rank tensor $\mathfrak{N}^{(s)}(r)$ are defined by formulas (\ref{eq5.18}), (\ref{eq5.31}) and (\ref{eq5.27}) respectively.
\end{lemma}

\subsection{Some functional spaces}

As a functional space $\mathbb{L}_{2}(\Omega)$ we denote the space of all measurable functions
$u(\boldsymbol{x})$ with a finite norm
\begin{equation*}
\|u\|_{2,\Omega}=\big(\int_{\Omega}u^{2}(\boldsymbol{x})dx\big)^{\frac{1}{2}}.
\end{equation*}
The functional space $\mathbb{W}^{1}_{2}(\Omega)$ is a closure of all infinitely smooth in $\Omega$ functions in the norm
\begin{equation*}
\|u\|^{(1)}_{2,\Omega}=\|u\|_{2,\Omega}+\|\nabla\,u\|_{2,\Omega}.
\end{equation*}
The functional space $\mathbb{C}^{k}(\overline{\Omega})$ for integer $k\geqslant 0$ consists of all
functions $u(\boldsymbol{x})$ with a finite norm
\begin{equation*}
|u|^{(k)}_{\Omega}=\max_{\boldsymbol{x}\in \Omega}\sum_{|m|=0}^{k}|D^{m}u(\boldsymbol{x})|,\,\,
D^{m}u=\frac{\partial^{|m|}u}{\partial x_{1}^{m_{1}}...\partial x_{n}^{m_{n}}},
\end{equation*}
where $m=(m_{1},...,m_{n}),\,m_{i}\geqslant 0,\,i=1,...,n$, is a multi-index, \\
$|m|=m_{1}+...+m_{n}$.

$\stackrel{\,\!\!\circ}{\mathbb{C}}^{\,k}(\overline{\Omega})$ is a subspace of
$\mathbb{C}^{k}(\overline{\Omega})$ of all functions $u(\boldsymbol{x})$ vanishing at the boundary $S=\partial\Omega$.

As a space $\mathbb{C}^{\infty}(\overline{\Omega})$ we denote the space of all infinitely smooth in $\Omega$ functions $u(\boldsymbol{x})$.

As a space $\mathbb{H}^{\alpha}(\overline{\Omega})$ we define the space of all functions $u(\boldsymbol{x})$  with a finite norm
\begin{equation*}
|u|^{(\alpha)}_{\Omega}=|u|^{(0)}_{\Omega_{T}}+\max_{\boldsymbol{x}\in\Omega,\,|h|<\infty}
\frac{|u(\boldsymbol{x}+h)-u(\boldsymbol{x})|}{|h|^{\alpha}}<\infty,\,\,
|u|^{(0)}_{\Omega_{T}}=\max_{(\boldsymbol{x},t)\in\Omega_{T}}|u(\boldsymbol{x},t)|.
\end{equation*}
As a space $\mathbb{H}^{\alpha,\frac{\alpha}{2}}(\overline{\Omega}_{t_{0}})$
we define the space of all functions $u(\boldsymbol{x},t)$ with a finite norm
\begin{equation*}
|u|^{(\alpha,\frac{\alpha}{2})}_{\Omega_{T}}=|u|^{(0)}_{\Omega_{T}}+
\max_{(\boldsymbol{x},t)\in\Omega_{T},\,|h|<\infty}
\frac{|u(\boldsymbol{x}+h,t+\frac{h}{2})-u(\boldsymbol{x},t)|}{|h|^{\alpha}}<\infty.
\end{equation*}
Finally, as a space $\mathbb{H}^{k+\alpha,\frac{k+\alpha}{2}}(\overline{\Omega}_{t_{0}})$ we define the space of all functions $u(\boldsymbol{x},t)$ with a finite norm
\begin{equation*}
|u|^{(k+\alpha,\frac{k+\alpha}{2})}_{\Omega_{T}}=|u|^{(0)}_{\Omega_{T}}+
\max_{(\boldsymbol{x},t)\in\Omega_{T},\,|h|<\infty}
\frac{|D^{k}u(\boldsymbol{x}+h,t+\frac{h}{2})-D^{k}u(\boldsymbol{x},t)|}{|h|^{\alpha}}.
\max_{(\boldsymbol{x},t)\in\Omega_{T}}<\infty.
\end{equation*}
As a space $\stackrel{\!\!\circ}{\mathbb{W}}^{1,0}_{2}(Q_{T})$ we define the space of all functions from
$\mathbb{W}^{1,0}_{2}(Q_{T}$) vanishing at the boundary $\partial{Q}$.

Let $\mathbb{W}$ be some functional space with elements $\boldsymbol{w}(\boldsymbol{x},t)$.
Then as $\mathbb{L}^{\infty}\big(0,T;\mathbb{W}(\Omega)\big)$ we denote all functions bounded in $\mathbb{W}(\Omega)$.

In particular, the space $\mathbb{L}^{\infty}\big((0,T;\mathbb{H}^{2+\alpha}\overline{\Omega})\big)$ consists of all bounded in $\mathbb{H}^{2+\alpha}(\overline{\Omega})$.

\section{Main Results}

\begin{theorem}\label{t3.1}

Let ${c}_{0}\in\mathbb{H}^{2+\alpha}(\overline{\Omega})$, $0\,\leqslant\,{c}_{0}\,\leqslant\,1$ and ${p}^{0}\in\mathbb{H}^{1+\alpha}(\overline{\Omega})$,
$\frac{\partial{p^{0}}}{\partial{x}_{3}}
(\boldsymbol{x})\,\geqslant\,0,\,\boldsymbol{x}\in\Omega$.

Then the problem $\mathbb{B}^{\varepsilon}(r)$ has a unique weak solution
$\widetilde{\boldsymbol{w}}^{\varepsilon}_{f},\,\widetilde{\boldsymbol{w}}^{\varepsilon}_{s},\,
\widetilde{c}\in\mathbb{W}^{1,0}_{2}(\Omega_{T})$, $\widetilde{p}^{\,\varepsilon}_{s}\in\mathbb{L}_{2}(\Omega_{T})$.

Moreover, due to the smoothness results (Ref~\cite{29}., Chapter~3, Theorem~1.1, \S~1)
\begin{equation}\label{eq3.1}
\widetilde{\boldsymbol{w}}_{s}^{\,\varepsilon}\in
\mathbb{L}^{\infty}\big(0,T;\mathbb{H}^{2+\alpha}(\overline{\Omega})\big).
\end{equation}
 \end{theorem}

\begin{theorem}\label{t3.2}
Under the conditions of the Theorem~\ref{t3.1} the problem $\mathbb{H}(r)$ has a unique weak solution
$\boldsymbol{w}_{f},\,\boldsymbol{w}_{s},\,{c}\in\mathbb{W}^{1,0}_{2}(\Omega_{T})$  ${p}_{f},\,{p}_{s}\in\mathbb{W}^{1,0}_{2}(\Omega_{T})$.
\end{theorem}
For the definition of weak solutions to the problem $\mathbb{B}^{\varepsilon}(r)$, see Section 2.
\
\begin{theorem}\label{t3.3}

Under the conditions of the Theorem~\ref{t3.1} the problem $\mathbb{H}$ has a unique solution
\begin{equation*}
{p}_{f}(.,t),\,{p}_{s},(.,t)\,\boldsymbol{w}_{f}(.,t),\,\boldsymbol{w}_{s}(.,t),\,
\in\mathbb{L}^{\infty}(0,T;\mathbb{H}^{2+\alpha}\overline{\Omega}),\,\,
{c}\in\mathbb{H}^{2+\alpha,\frac{2+\alpha}{2}}(\overline{\Omega}_{T})
\end{equation*}
and
\begin{equation}\label{eq3.2}
\max_{t\in[0,T]}\big(|p_{f}(.,t)|_{\Omega}^{(2+\alpha)},\,|p_{s}(.,t)|_{\Omega}^{(2+\alpha)},\,
|\boldsymbol{w}_{f}(.,t)|_{\Omega}^{(2+\alpha)},\,|\boldsymbol{w}_{s}(.,t)|_{\Omega}^{(2+\alpha)}\big)
\leqslant\,M^{(0)}
\end{equation}
\begin{equation}\label{eq3.3}
|{c}|_{\Omega_{T}}^{(2+\alpha)}\,\leqslant\,M^{(0)},
\end{equation}
where $M^{(0)}$ depends on given data ${c}^{0}(\boldsymbol{x})$ and ${p}^{0}(\boldsymbol{x})$.
\end{theorem}

\section{Proof of Theorem~\ref{t3.1}: the Existence of a Weak Solution to the Problem $\mathbb{B}^{\varepsilon}(r)$}

\subsection{Proof of Theorem~\ref{t3.1}: existence of a weak solution to the dynamic problem $\mathbb{B}^{\varepsilon}(r)$}

Thanks to the linearity of the problem $\mathbb{B}^{\varepsilon}(r)$, it is sufficient to derive
the corresponding a priori estimates.

\begin{lemma} \label{l4.1}
Under the conditions of Theorem~\ref{t2.1} the dynamic problem $\mathbb{B}^{\varepsilon}(r)$ has a unique
weak solution such that
\begin{multline}\label{eq4.1}
\max_{0<t<T}\big(\|\chi^{\,\varepsilon}(.,t)
\big(\widetilde{\boldsymbol{w}}^{\,\varepsilon}_{f}(.,t)\|_{2,\Omega}+
\|\sqrt{\varepsilon}\frac{\partial\widetilde{\boldsymbol{w}}^{\,\varepsilon}_{f}}{\partial{t}}
\|_{2,\Omega_{T}}+
\\
\max_{0<t<T}\varepsilon\|\chi^{\varepsilon}(.,t)
\mathbb{D}\big(x,\widetilde{\boldsymbol{w}}^{\,\varepsilon}_{f}(.,t)\big)\|_{2,\Omega}+
\|\mathbb{D}(x,\varepsilon
\frac{\partial\widetilde{\boldsymbol{w}}^{\,\varepsilon}_{f}}{\partial{t}})\|_{2,\Omega_{T}}
\leqslant\,{M},
\end{multline}
\begin{equation}\label{eq4.2}
\max_{0<t<T}\Big(\|\big(1-\chi^{\varepsilon}(.,t)
\widetilde{\boldsymbol{w}}^{\,\varepsilon}_{s}(.,t)\big)\|_{2,\Omega}+
\|(1-\chi^{\varepsilon}(.,t)\mathbb{D}\big(x,\widetilde{\boldsymbol{w}}^{\,\varepsilon}_{s}
(.,t)\big)\|_{2,\Omega}\leqslant\,{M},
\end{equation}
\begin{equation}\label{eq4.3}
\max_{0<t<T}(\|\chi^{\,\varepsilon}
\nabla\cdot\widetilde{\boldsymbol{w}}^{\,\varepsilon}_{f}(.,t_{0})\|_{2,\Omega}+
\|(1-\chi^{\varepsilon})\nabla\cdot\widetilde{\boldsymbol{w}}^{\,\varepsilon}_{s}(.,t_{0})\|_{2,\Omega}
\leqslant\,{M},
\end{equation}
where $M$ does not depend on $\varepsilon$.
\end{lemma}

\begin{proof}
Let in (\ref{eq2.43}) $\boldsymbol{\varphi}=(1-\chi^{\varepsilon})\widetilde{\boldsymbol{w}}^{\,\varepsilon}_{s}$.

Then using the continuity equation (\ref{eq2.4}), the simplest embedding theorem (Lemma~\ref{l2.2})
and {H\"older}'s inequality (Lemma~\ref{l2.5}) we obtain
\begin{multline*}
\lambda_{0}\int_{0}^{t_{0}}\int_{\Omega}\big((1-\chi^{\varepsilon}\big)
|\mathbb{D}\big(x,\widetilde{\boldsymbol{w}}^{\,\varepsilon}_{s}\big)|^{2}dxdt+
\int_{0}^{t_{0}}\int_{\Omega}\big((1-\chi^{\varepsilon})\big){c}^{2}_{s}\varrho^{0}_{s}
|\nabla\cdot\widetilde{\boldsymbol{w}}^{\,\varepsilon}_{s}|^{2}\big)dx=
\\
|\int_{0}^{t_{0}}\int_{\Omega}(1-\chi^{\varepsilon})
(\nabla{p}_{0}\cdot\widetilde{\boldsymbol{w}}^{\,\varepsilon}_{s})dxdt|\leqslant
\\
\frac{\delta}{2}\int_{\Omega}\big((1-\chi^{\varepsilon})\big)
|\widetilde{\boldsymbol{w}}^{\,\varepsilon}_{s}(.,)|^{2}dx+
\frac{2}{\delta}\int_{\Omega}|\nabla{p}_{0}|^{2}dx\leqslant
\\
\frac{\delta}{2}\int_{\Omega}\big((1-\chi^{\varepsilon})(.,{t}_{0})\big)
|\mathbb{D}\big(x,\widetilde{\boldsymbol{w}}^{\,\varepsilon}_{s}(.,{t}_{0})|^{2}dx+\frac{{M}}{\delta}.
\end{multline*}
To estimate liquid displacements we put in (\ref{eq2.48})
$\displaystyle\,\frac{\partial\boldsymbol{\varphi}}{\partial{t}}=\chi^{\varepsilon}
\frac{\partial\widetilde{\boldsymbol{w}}^{\,\varepsilon}_{f}}{\partial{t}}$ and using continuity equation (\ref{eq2.2}), trivial inequality 
$\displaystyle\,{a}{b}\leqslant\frac{\delta}{2}{a}+\frac{1}{2\delta}{b}^{2}$ and integration by parts 
we arrive at
\begin{multline*}
\varepsilon\int_{0}^{t_{0}}\int_{\Omega}\chi^{\varepsilon}
|\frac{\partial\widetilde{\boldsymbol{w}}^{\,\varepsilon}_{f}}{\partial{t}}|^{2}dxdt+
\int_{\Omega}\chi^{\varepsilon}\varepsilon^{2}\mu_{1}|\mathbb{D}
\big(x,\frac{\partial\widetilde{\boldsymbol{w}}^{\varepsilon}_{f}}{\partial{t}}(.,t_{0})\big)|^{2}dx=
\\
\int_{0}^{t_{0}}\int_{\Omega}\chi^{\varepsilon}\varepsilon
|\frac{\partial\widetilde{\boldsymbol{w}}^{\,\varepsilon}_{f}}{\partial{t}}|^{2}dxdt+
\int_{\Omega}\chi^{\varepsilon}(.,t_{0})\varepsilon^{2}\mu_{1}
|\mathbb{D}\big(x,\widetilde{\boldsymbol{w}}^{\,\varepsilon}_{f}(.,t_{0})\big)|^{2}dx
\\
=|\int_{\Omega}\chi^{\varepsilon}(.,t_{0})
\big(\widetilde{\boldsymbol{w}}^{\,\varepsilon}_{f}(.,t_{0})\cdot\nabla{p}_{0}\big)dx|\leqslant\,
\\
\frac{\delta}{2}\int_{\Omega}\chi^{\varepsilon}(.,t_{0})
|\widetilde{\boldsymbol{w}}^{\,\varepsilon}_{f}(.,t_{0})|^{2}dx+
\frac{1}{2\delta}\int_{\Omega}|\nabla{p}^{0}|^{2}dx.
\end{multline*}
Next we apply the Poincar\'{e} inequality (\ref{eq2.29})
\begin{equation*}
\int_{\Omega}\chi^{\varepsilon}(.,t)|\widetilde{\boldsymbol{w}}^{\,\varepsilon}_{f}(.,t)|^{2}dx\,\leqslant
{M}_{\Omega}\mu_{1}^{-1}\int_{\Omega}\chi^{\varepsilon}(.,t)\varepsilon^{2}\mu_{1}
\mathbb{D}\big(x,\widetilde{\boldsymbol{w}}^{\,\varepsilon}_{f}(.,t)\big)|^{2}dx
\end{equation*}
and obtain
\begin{multline*}
\int_{0}^{t_{0}}\int_{\Omega}\chi^{\varepsilon}\varepsilon
|\frac{\partial\widetilde{\boldsymbol{w}}^{\,\varepsilon}_{f}}{\partial{t}}|^{2}dxdt+
\int_{0}^{t_{0}}\int_{\Omega}\chi^{\varepsilon}\varepsilon^{2}\mu_{1}
|\mathbb{D}(x,\frac{\partial\widetilde{\boldsymbol{w}}^{\,\varepsilon}_{f}}{\partial{t}})|^{2}dxdt+
\\
\int_{\Omega}\chi^{\varepsilon}(.,t_{0})\varepsilon^{2}\mu_{1}
|\mathbb{D}\big(x,\widetilde{\boldsymbol{w}}^{\,\varepsilon}_{f}(.,t_{0})\big)|^{2}dx+
\int_{\Omega}\chi^{\,\varepsilon}|\nabla\cdot
\widetilde{\boldsymbol{w}}_{f}^{\,\varepsilon}(.,t_{0})|^{2}dx\leqslant
\\
\frac{\delta}{2}{M}_{\Omega}\mu_{1}^{-1}\int_{\Omega}\chi^{\varepsilon}\varepsilon^{2}\mu_{1}
|\mathbb{D}(x,\frac{\partial\widetilde{\boldsymbol{w}}^{\,\varepsilon}_{f}}{\partial{t}})|^{2}dx+
\frac{1}{2\delta}\int_{\Omega}|\nabla{p}^{0}|^{2}dx.
\end{multline*}
The desired estimates follow from the last inequality for $\delta=2(\mu_{1}{M}_{\Omega})^{-1}$.
\end{proof}

\subsection{Proof of Theorem~\ref{t3.1}: the existence of a weak solution to the diffusion problem $\mathbb{B}^{\varepsilon}(r)$}
\begin{lemma} \label{l4.2}
Under the conditions of Theorem~\ref{t2.1} the diffusion problem $\mathbb{B}^{\varepsilon}(r)$ has a unique weak solutions $\widetilde{{c}}^{\,\varepsilon}$, such that
\begin{equation}\label{eq4.4}
\|(\widetilde{{c}}^{\,\varepsilon}-{c}^{0})\|_{2,\Omega_{T}}+
\|\nabla(\widetilde{{c}}^{\,\varepsilon}-{c}^{0})\|_{2,\Omega_{T}}
\leqslant\,{M}\|\nabla{c}^{0}\|_{2,\Omega},
\end{equation}
where $M$ does not depend on $\varepsilon$.
\end{lemma}
\begin{proof}
To prove this, we only need to obtain a priori estimates to the solution of the diffusion problem $\mathbb{B}^{\varepsilon}(r)$, written in the equivalent form of the integral identity (\ref{eq2.45}). 
To do this we repeat the proof of Lemma~2.1 in \S~2, chapter~III of~\cite{28} with test function $\xi=\widetilde{{c}}^{\,\varepsilon}-{c}^{0}$ using trivial inequality $|ab|\leqslant\,\delta{a}^{2}+\frac{b^{2}}{4\delta}$ for any $\delta>0$, the H\"older's inequality (\ref{eq2.35}) and integration by parts we obtain the chain of inequalities%

\begin{multline*} 
0=\int_{\Omega}\chi^{\varepsilon}(.,t_{0})\big(\widetilde{{c}}^{\,\varepsilon}
(\boldsymbol{x},t_{0})-{c}^{0}(\boldsymbol{x})+
\frac{\beta^{\varepsilon}}{\alpha^{\varepsilon}}+{c}^{0}(\boldsymbol{x}))
(\widetilde{{c}}^{\,\varepsilon}(\boldsymbol{x},t_{0})-
{c}^{0}(\boldsymbol{x})\big)dx-
\\
\int_{0}^{t_{0}}\int_{\Omega}\chi^{\varepsilon}(.,t_{0})
(\widetilde{{c}}^{\,\varepsilon}-{c}^{0}+
\frac{\beta^{\varepsilon}}{\alpha^{\varepsilon}}+{c}^{0})
\frac{\partial}{\partial{t}}(\widetilde{{c}}^{\,\varepsilon}-{c}^{0}+{c}^{0}+
\frac{\beta^{\varepsilon}}{\alpha^{\varepsilon}})dxdt+
\\
\alpha_{0}\int_{0}^{t_{0}}\int_{\Omega}\chi^{\varepsilon}
\big(\nabla(\widetilde{{c}}^{\,\varepsilon}-{c}^{0})
\cdot\nabla(\widetilde{{c}}^{\,\varepsilon}-{c}^{0}+{c}^{0})\big)dxdt=
\\
\int_{\Omega}\chi^{\varepsilon}(.,t_{0})
\Big(\big(\widetilde{{c}}^{\,\varepsilon}(.,t_{0})-{c}^{0}\big)^{2}+
(\frac{\beta^{\varepsilon}}{\alpha^{\varepsilon}}+{c}^{0})
\big(\widetilde{{c}}^{\,\varepsilon}(.,t_{0})-{c}^{0}\big)\Big)dx-
\\
-\frac{1}{2}\int_{0}^{t_{0}}\int_{\Omega}\chi^{\varepsilon}
\frac{\partial}{\partial{t}}(\widetilde{{c}}^{\,\varepsilon}(.,t_{0})-{c}^{0}+{c}^{0}+
\frac{\beta^{\varepsilon}}{\alpha^{\varepsilon}})^{2}dxdt
+\alpha_{c}\int_{0}^{t_{0}}\int_{\Omega^{\varepsilon}_{f}(r(.,t))}
|\nabla(\widetilde{{c}}^{\,\varepsilon}-{c}_{0})|^{2}dxdt+
\\
\alpha_{c}\int_{0}^{t_{0}}\int_{\Omega}\chi^{\varepsilon}
\big(\nabla(\widetilde{{c}}^{\,\varepsilon}-{c}^{0})\cdot\nabla{c}^{0}\big)dxdt=
\\
\int_{\Omega}\chi^{\varepsilon}(.,t_{0})
\Big(\big(\widetilde{{c}}^{\,\varepsilon}(.,t_{0})-{c}^{0}\big)^{2}+
(\frac{\beta^{\varepsilon}}{\alpha^{\varepsilon}}+{c}^{0})
\big(\widetilde{{c}}^{\,\varepsilon}(.,t_{0})-{c}^{0})+
\frac{1}{2}(\widetilde{{c}}_{0}+
\frac{\beta^{\varepsilon}}{\alpha^{\varepsilon}})^{2}\big)\Big)dx+
\\
\frac{1}{2}\int_{0}^{t_{0}}\int_{\Gamma^{\varepsilon}(r(.,t))}
(\widetilde{{c}}^{\,\varepsilon}+
\frac{\beta^{\varepsilon}}{\alpha^{\varepsilon}})^{2}D^{\varepsilon}_{N}
\sin\psi\,d\sigma{d}t+
\alpha_{c}\int_{0}^{t_{0}}\int_{\Omega^{\varepsilon}_{f}(r(.,t))}
|\nabla(\widetilde{{c}}^{\,\varepsilon}-{c}^{0}|^{2}dxdt+
\\
\alpha_{c}\int_{0}^{t_{0}}\int_{\Omega}\big(\nabla(\widetilde{{c}}^{\,\varepsilon}-{c}^{0})\cdot
\nabla{c}^{0}\big)dxdt\geqslant
\\
\int_{\Omega}\chi^{\varepsilon}(.,t_{0})
\Big(\big(\widetilde{{c}}^{\,\varepsilon}(.,t_{0})-{c}^{0}\big)^{2}+
(\frac{\beta^{\varepsilon}}{\alpha^{\varepsilon}})
\big(\widetilde{{c}}^{\,\varepsilon}(.,t_{0})-{c}^{0})+\frac{1}{2}({c}^{0}+
\frac{\beta^{\varepsilon}}{\alpha^{\varepsilon}})^{2}\big)\Big)dx+
\\
\alpha_{c}\int_{0}^{t_{0}}\int_{\Omega}\chi^{\varepsilon}
|\nabla(\widetilde{{c}}^{\,\varepsilon}-{c}^{0})|^{2}dxdt+
\alpha_{c}\int_{0}^{t_{0}}\int_{\Omega^{\varepsilon}_{f}(r(.,t))}
\big(\nabla(\widetilde{{c}}^{\,\varepsilon}-{c}^{0})\cdot
\nabla{c}^{0}\big)dxdt\geqslant
\end{multline*}
\begin{multline}\label{eq4.5} 
\\
\int_{\Omega}\chi^{\varepsilon}(.,t_{0})
\Big(\big(\widetilde{{c}}^{\,\varepsilon}(.,t_{0})-{c}^{0}\big)^{2}dxdt+
\frac{\alpha_{c}}{2}\int_{0}^{t_{0}}\int_{\Omega^{\varepsilon}_{f}(r(.,t))}
|\nabla(\widetilde{{c}}^{\,\varepsilon}-{c}^{0})|^{2}dxdt-
\\
\frac{\alpha_{c}}{2}\int_{\Omega}|\nabla{c}^{0}|^{2}dxdt
\end{multline}%

which proves the statement of the lemma.
\end{proof}

\section{Proof of Theorem~\ref{t3.2}: Homogenization of the Problem $\mathbb{B}^{\varepsilon}(r)$}

The homogenization procedure itself is well explained in many publications (\cite{6,7,8,9,10,11,12,13,16,36}). For the dynamic problem, the reader can follow the proof of Theorem~1 in chapter~I, section~1.3 of~\cite{16} and for the diffusion problem, chapter~10 of~\cite{16}.

\begin{lemma} \label{l5.1}
Under the conditions of  the Lemma~\ref{l4.1} there exist functions $\boldsymbol{w}_{f}$,
$\boldsymbol{w}_{s}$, ${p}$, $\pi$, ${c}$ and 1-periodic in the variable $\boldsymbol{y}$ functions $\boldsymbol{W}_{f}(\boldsymbol{y};\boldsymbol{x},t)$, $\mathbb{D}\big(y,\boldsymbol{W}_{f}(\boldsymbol{y};\boldsymbol{x},t)\big)$, $\nabla\cdot\boldsymbol{w}_{f}$,
$\boldsymbol{W}_{s}(\boldsymbol{y};\boldsymbol{x},t)$,   $\mathbb{D}\big(y,\boldsymbol{W}_{s}(\boldsymbol{y};\boldsymbol{x},t)\big)$, $\nabla\cdot\boldsymbol{w}_{s}$,
$\Pi(\boldsymbol{y};\boldsymbol{x},t)$ and $C(\boldsymbol{y};\boldsymbol{x},t)$ such that
$\boldsymbol{w}_{f}\in\mathbb{L}_{2}(\Omega_{T})$, $\boldsymbol{w}_{s},\,{\pi},\,{c}\,\in\mathbb{W}^{1,0}_{2}(\Omega_{T})$,
$\boldsymbol{W}_{f},\,{C}\in\mathbb{L}_{2}\big(0,T;\mathbb{W}^{1}_{2}(\textbf{Y}_{f})\big)$
and $\boldsymbol{W}_{s}\in\mathbb{L}_{2}(\textbf{Y}_{s})$.

1) The sequence $\{\widetilde{\boldsymbol{w}}^{\,\varepsilon}_{f}\}$ converges weakly to the function $\boldsymbol{w}_{f}$ and two-scale to the function $\boldsymbol{W}_{f}(\boldsymbol{y};\boldsymbol{x},t)$.

2) The sequences $\{\varepsilon\mathbb{D}(x,\widetilde{\boldsymbol{w}}^{\,\varepsilon}_{f})\}$ and
$\{\varepsilon\nabla_{x}\cdot\widetilde{\boldsymbol{w}}^{\,\varepsilon}_{f})\}$ converge two-scale to the functions $\mathbb{D}(y,\boldsymbol{W}_{f})$ and $\nabla_{y}\cdot\boldsymbol{W}_{f}$ respectively.

3) The sequence $\{\widetilde{\boldsymbol{w}}^{\,\varepsilon}_{s}\}$, converge two-scale and
weakly to the function $\boldsymbol{w}_{s}\in\mathbb{L}_{2}(\Omega_{T})$.

4) The sequences $\{\mathbb{D}(x,\widetilde{\boldsymbol{w}}^{\,\varepsilon}_{s})\}$ and $\{\nabla\cdot\widetilde{\boldsymbol{w}}^{\,\varepsilon}_{s}\}$ converge two-scale to the functions $\mathbb{D}(x,\boldsymbol{w}_{s})+D(y,\boldsymbol{W}_{s})$ and $\nabla\cdot\boldsymbol{w}_{s}+\nabla_{y}\cdot\boldsymbol{W}_{s}$.

5) The sequence $\{\widetilde{c}^{\,\varepsilon}\}$ converges weakly and two-scale to the function
${c}\in\mathbb{W}^{1,0}_{2}(\Omega_{T})$.

6) The sequence $\{\nabla\widetilde{c}^{\,\varepsilon}\}$ converges two-scale to the function
${\nabla{c}}+\nabla_{y}{C}$.

Here $\boldsymbol{W}_{s}\in\mathbb{L}_{2}(\Omega_{T})\cap\mathbb{W}^{1,0}_{2}(\textbf{Y}_{s})$,
$\boldsymbol{W}_{f}$, ${C}$, ${\Pi}\in\mathbb{L}_{2}(\Omega_{T})\cap\mathbb{W}^{1,0}_{2}(\textbf{Y}_{f})$.

7) Hold true the following a priori estimates
\begin{multline}\label{eq5.1}
\,\,\,\,\,\,\,\,\,\,\,\,\,\,\,\,\,\,\,\,\,\,\,\,\,\,\,\,\,\,
\|{\boldsymbol{w}}_{f}\|_{2,\Omega_{T}}+\|\boldsymbol{w}_{s}\|^{(1,0)}_{2,\Omega_{T}}+
\|\boldsymbol{W}_{f}\|_{2,\textbf{Y}_{f}\times\Omega_{T}}+
\\
\|\mathbb{D}\big(y,\boldsymbol{W}_{s})\|_{2,{\textbf{Y}}_{s}\times\Omega_{T}}+
\|\mathbb{D}\big(y,\boldsymbol{W}_{f})\|_{2,{Y}_{f}\times\Omega_{T}}\leqslant\,{M},
\,\,\,\,\,\,\,\,\,\,\,\,\,\,\,
\end{multline}
\begin{equation}\label{eq5.2}
\|({C}-{c}_{0})\|^{(1,0)}_{2,Y\times\Omega_{T}}+\|({c}-{c}_{0})\|^{(1,0)}_{2,\Omega_{T}},
\end{equation}
\begin{equation}\label{eq5.3}
\|{p}_{f}\|_{2,\Omega_{T}}+\|{p}_{s}\|_{2,\Omega_{T}}+
\|{\pi}\|_{2,\Omega_{T}}+|\frac{\partial{\pi}}{\partial{t}}\|_{2,\Omega_{T}}\leqslant\,{M},
\end{equation}
where $M$ does not depend on $\varepsilon$.
\end{lemma}
The proof is straightforward and based on the estimates (\ref{eq4.1})-(\ref{eq4.4}).

We only note that
\begin{equation*}
\varepsilon\int_{0}^{t_{0}}\int_{\Omega}\chi^{\varepsilon}
|\frac{\partial\widetilde{\boldsymbol{w}}^{\,\varepsilon}_{f}}{\partial{t}}|^{2}dxdt\leqslant\,{M}
\end{equation*}
and
\begin{equation*}
\lim_{\varepsilon\rightarrow{0}}\varepsilon\int_{0}^{t_{0}}\int_{\Omega}\chi^{\varepsilon}
|\frac{\partial\widetilde{\boldsymbol{w}}^{\,\varepsilon}_{f}}{\partial{t}}|dxdt=0.
\end{equation*}

\begin{remark}\label{5.1}
We have reconstructed the liquid pressure ${p}_{f}$ and its antiderivative ${\pi}_{f}$  by decomposing the space $\mathbb{L}_{2}(\Omega)$ into a direct sum of the subspace of all
solenoidal vector functions and the subspace of all gradients of scalar functions (see Lemma \ref{l2.5}).
The solid pressure ${p}_{s}$ is restored after the two-scale limit in the continuity equation (\ref{eq2.47}).
\end{remark}

\subsection{Proof of Theorem~\ref{t3.3}: homogenization of the dynamic problem $\mathbb{B}^{\varepsilon}(r)$ for the liquid component}

\begin{lemma} \label{l5.2}
Under the conditions of the Theorem~\ref{t2.1} the limiting procedure in the equations (\ref{eq2.1}) and (\ref{eq2.2}) and the integral identity (\ref{eq2.47}) results in the following dynamic problem $\mathbb{H}(r)$ for displacements and pressure of the liquid component consisting of Darcy law of filtration
\begin{equation}\label{eq5.4}
\boldsymbol{w}_{f}=\int_{\textbf{Y}_{f}}\boldsymbol{W}_{f}dy=
-\frac{1}{\mu_{1}}(B)^{(w)}(r)<\nabla_{x}{\pi}>,\,\,
\nabla_{x}\cdot\boldsymbol{w}_{f}=0,
\end{equation}
for the liquid displacements $\boldsymbol{w}_{f}$ and the antiderivative ${\pi}$  of the pressure ${p}$ in the domain $\Omega_{T}$, completed with the boundary conditions
\begin{equation}\label{eq5.5}
{\pi}(\boldsymbol{x},t)={p}^{0}t,\,\,\,\boldsymbol{x}\in{S}^{1}\cup{S}^{2},\,\,\,0<t<T,
\end{equation}
\begin{equation}\label{eq5.6}
\boldsymbol{w}_{f}\cdot\boldsymbol{n}=0,\,\,\,\boldsymbol{x}\in{S}^{0},\,\,\,0<t<T,
\end{equation}
where $\boldsymbol{n}$ is a normal vector to the boundary ${S}^{0}$.

The symmetric strictly positive definite constant matrix $\mathbb{B}^{(w)}(r)$ is given by formula (\ref{eq5.18}).
\end{lemma}

\begin{proof}
First, we derive the continuity equations for functions $\boldsymbol{w}_{f}$ and $\boldsymbol{W}_{f}$.

To do that we first consider continuity equation \eqref{eq2.2} in the form multiply it by function ${\eta}$, that vanishes in ${S}^{1}\cup{S}^{2}$, integrate by parts and take the limit as $\varepsilon\rightarrow{0}$:
\begin{equation*}
0=\lim_{\varepsilon\rightarrow{0}}\int_{0}^{{t}_{0}}\int_{\Omega}{\eta}
\widetilde{\boldsymbol{w}}_{f}\cdot\nabla{\eta}\big)dxdt.
\end{equation*}
This identity gives us
\begin{equation}\label{eq5.7}
\nabla\cdot\boldsymbol{w}_{f}=0.
\end{equation}
To derive the continuity equation for the unknown functions $\boldsymbol{W}_{f}(\boldsymbol{y};\boldsymbol{x},t)$ (liquid displacements) and $\boldsymbol{\Pi}^{(f)}(\boldsymbol{y};\boldsymbol{x},t)$ (liquid pressure) we consider the integral identity \eqref{eq2.2} with arbitrary test functions ${\xi}=\varepsilon{\eta}(\boldsymbol{x},t)\phi(\frac{\boldsymbol{x}}{\varepsilon})$, where ${\eta}(\boldsymbol{x},t)$ is the same as before and $\phi(\boldsymbol{y})$ is 1-periodic in $\boldsymbol{y}$ function, and using relations 1) and 2) in conditions of the Lemma~\ref{l5.1} obtain:
\begin{equation}\label{eq5.8}
0=\lim_{\varepsilon\rightarrow{0}}\int_{0}^{{t}_{0}}
\int_{\Omega}{\eta}\phi\chi^{\varepsilon}
\varepsilon\nabla\cdot\widetilde{\boldsymbol{w}}_{f}^{\varepsilon}dxdt=
\int_{0}^{{t}_{0}}\int_{\Omega}{\eta}\int_{Y_{f}(r)}
\big(\phi\nabla_{y}\cdot\boldsymbol{W}_{f}\big)dydxdt.
\end{equation}
Due to the arbitrary choice of the functions ${\eta}$ and $\phi$ the last relation means that holds true the continuity equation
\begin{equation}\label{eq5.9}
\nabla_{y}\cdot\boldsymbol{W}_{f}(\boldsymbol{y};\boldsymbol{x},t)=0,\,\,\,
(\boldsymbol{y};\boldsymbol{x},t)\in{Y}_{f}\times\Omega_{t_{0}}
\end{equation}
and boundary and normalization conditions
\begin{equation}\label{eq5.10}
(\boldsymbol{W}\cdot\boldsymbol{N})=0,\,\,\int_{\textbf{Y}_{f}}\boldsymbol{W}dy=0,\,\,
(\boldsymbol{y};\boldsymbol{x},t)\in\gamma\times\Omega_{t_{0}},\,\,j=ol,\,sp.
\end{equation}
Let additionally in (\ref{eq2.44}) $(1-\chi^{\varepsilon})\boldsymbol{\varphi}=0$ and
$\boldsymbol{\varphi}=0$ at $t=0$ and $t=t_{0}$,  $\frac{\partial\boldsymbol{\varphi}}{\partial{t}}=
\eta(\boldsymbol{x},t)\boldsymbol{\psi}(\frac{\boldsymbol{x}}{\varepsilon})$, where
$\eta\in{W}^{1,1}_{2}(\Omega_{T})$, $\eta(\boldsymbol{x},t)=0$ for $\boldsymbol{x}\in{S}^{0}$,
$0<t<T$ and $\boldsymbol{\psi}\in{W}^{1}_{2}(Y_{f})$,
\emph{supp}$\boldsymbol{\psi}\subset{Y}_{f}$, $\nabla_{y}\cdot\boldsymbol{\psi}=0$.

Then
\begin{multline*}
\mathbb{D}(x,{\eta}\boldsymbol{\psi})=\frac{1}{2}\sum_{i,j=1}^{3}\big({d}_{ij}(x,{\eta}\boldsymbol{\psi})
\boldsymbol{e}^{i}\otimes\boldsymbol{e}^{j}+{d}_{ji}(x,{\eta}\boldsymbol{\psi})
\boldsymbol{e}^{j}\otimes\boldsymbol{e}^{i}\big)=
\\
\frac{{\eta}}{2\varepsilon}\sum_{i,j=1}^{3}(\frac{\partial{\psi}_{i}}{\partial{y}_{j}}
\boldsymbol{e}^{i}\otimes\boldsymbol{e}^{j}+\frac{\partial{\psi}_{j}}{\partial{y}_{i}}
\boldsymbol{e}^{j}\otimes\boldsymbol{e}^{i}))+\frac{1}{2}\sum_{i,j=1}^{3}({\psi}_{i}
\frac{\partial{\eta}}{\partial{x}_{j}}\boldsymbol{e}^{i}\otimes\boldsymbol{e}^{j}+{\psi}_{j}
\frac{\partial{\eta}}{\partial{x}_{i}}\boldsymbol{e}^{j}\otimes\boldsymbol{e}^{i}),
\\
\varepsilon^{2}\mathbb{D}(x,\frac{\partial\boldsymbol{\varphi}}{\partial{t}})=
\eta\varepsilon\mathbb{D}\big(y,\boldsymbol{\psi}(\frac{\boldsymbol{x}}{\varepsilon})\big)+
\frac{\varepsilon^{2}}{2}(\nabla{\eta}\otimes\boldsymbol{\psi}+
\boldsymbol{\psi}\otimes\nabla{\eta}),\,\,\,\nabla\cdot(\eta\boldsymbol{\psi})=
(\nabla{\eta}\cdot\boldsymbol{\psi})=-\frac{1}{\varepsilon}\eta(\nabla_{y}\cdot\boldsymbol{\psi})=0.
\end{multline*}
Next we consider functions
\begin{multline}\label{eq5.11}
A_{f}(\boldsymbol{x},t)=\int_{Y_{f}}\big(\nabla({p}^{0}t))\cdot\boldsymbol{\psi}\big)-\nabla_{y}
\cdot\big(\mu_{1}\mathbb{D}(y,\boldsymbol{W}_{f})\big)dy=
-\int_{Y_{f}}\nabla_{y}\cdot\big(\mu_{1}\mathbb{D}(y,\boldsymbol{W}_{f})\big)dy\,\,
\\
({B})_{f}(\boldsymbol{x},t)=(\int_{Y_{f}}\boldsymbol{\psi}dy)<t{\pi}(\boldsymbol{x},t)-{p}^{0}>
\end{multline}
and functional
\begin{equation*}
I^{\varepsilon}_{f}(\eta\boldsymbol{\psi})=\int_{0}^{t_{0}}\int_{\Omega}\chi^{\varepsilon}\eta\int_{Y_{f}}
\Big(\big(\varepsilon\widetilde{\boldsymbol{w}}_{f}^{\varepsilon})\cdot\boldsymbol{\psi}\big)+
\big(\mu_{1}\varepsilon^{2}\mathbb{D}(x,\widetilde{\boldsymbol{w}}_{f}^{\varepsilon})-
\\
(t\widetilde{{\pi}}^{\varepsilon}-{p}^{0})\mathbb{I}\big):\mathbb{D}(x,\boldsymbol{\psi})\Big)dydxdt.
\end{equation*}
In accordance with Lemma~\ref{l5.1} and the integral identity (\ref{eq2.44}) we get
\begin{multline}\label{eq5.12}
0={I}^{0}_{f}(\eta\boldsymbol{\psi})=
\lim_{\varepsilon\rightarrow{0}}I^{\varepsilon}_{f}(\eta\boldsymbol{\psi})=
-\lim_{\varepsilon\rightarrow{0}}\int_{0}^{t_{0}}\int_{\Omega}\chi^{\varepsilon}
\Big((\varepsilon\widetilde{\boldsymbol{w}}_{f}^{\varepsilon}\cdot
\frac{\partial\boldsymbol{\varphi}}{\partial{t}})+
\\
\big(\mu_{1}\varepsilon^{2}\mathbb{D}(x,\widetilde{\boldsymbol{w}}_{f}^{\varepsilon})+
(t\widetilde{{\pi}}^{\varepsilon}-{p}^{0})\mathbb{I}\big):
\mathbb{D}(x,\frac{\partial\boldsymbol{\varphi}}{\partial{t}})\Big)dxdt=
\\
-\lim_{\varepsilon\rightarrow{0}}\int_{0}^{t_{0}}\int_{\Omega}\eta\chi^{\varepsilon}
\Big(\big(\nabla(\varepsilon\widetilde{\boldsymbol{w}}_{f}^{\varepsilon})
\cdot\boldsymbol{\psi}\big)+\mu_{1}\varepsilon\mathbb{D}(x,\widetilde{\boldsymbol{w}}_{f}^{\varepsilon}):
\mathbb{D}\big(y,\boldsymbol{\psi}(\frac{\boldsymbol{x}}{\varepsilon})\big)-
(t\widetilde{{\pi}}^{\varepsilon}-{p}^{0})\big(\nabla\eta\cdot
\boldsymbol{\psi}(\frac{\boldsymbol{x}}{\varepsilon})\big)\Big)dxdt=
\\
-\int_{0}^{t_{0}}\int_{\Omega}\Big(\eta\big(\int_{Y_{f}}\big(
\mu_{1}\mathbb{D}(y,\boldsymbol{W}_{f}):\mathbb{D}(y,\boldsymbol{\psi})dy+
(t{\pi}-{p}^{0})(\int_{Y_{f}}\boldsymbol{\psi}dy\cdot\nabla\eta)\Big)dxdt=
\\
\int_{0}^{t_{0}}\int_{\Omega}\big(({B})_{f}<\nabla\eta>-A_{f}\eta\big)dxdt=0.
\end{multline}
The last identity in (\ref{eq5.12})
\begin{equation}\label{eq5.13}
\int_{0}^{t_{0}}\int_{\Omega}\big(({B})_{f}<\nabla\eta>-A_{f}\eta\big)dxdt=0
\end{equation}
means that function ${\pi}\in\mathbb{W}^{1,0}_{2}(\Omega_{T})$ and identity
(\ref{eq5.13}) takes the form of the differential equation
\begin{equation}\label{eq5.14}
\nabla_{y}\cdot\big(\mu_{1}\mathbb{D}(y,\boldsymbol{W}_{f})\big)=
-\nabla_{x}(t\pi-{p}^{0})(\boldsymbol{x},t)=
-\sum_{i=1}^{3}\frac{\partial{}}{\partial{x_{i}}}({\pi}-{p}^{0}t)(\boldsymbol{x},t)\boldsymbol{e}^{i},
\end{equation}
completed with the continuity equation (\ref{eq5.9}), boundary condition (\ref{eq5.11}) and boundary condition (\ref{eq5.5})
\begin{equation}\label{eq5.15}
t\pi(\boldsymbol{x},t)-{p}^{0}=0,\,\,\,\boldsymbol{x}\in{S}^{1}\cup{S}^{2},\,\,0<t<T,
\end{equation}
which is a consequence of the identity (\ref{eq5.13}).

To solve the periodic boundary value problem (\ref{eq5.14}), (\ref{eq5.9}) and (\ref{eq5.15}) we use decomposition
\begin{multline}\label{eq5.16}
\,\,\,\,\,\,\,\,\,
\boldsymbol{w}_{f}=\int_{\textbf{Y}_{f}}\boldsymbol{W}_{f}\big(\boldsymbol{y};\boldsymbol{x},t\big)dy=
\\
-\frac{1}{\mu_{1}}\sum_{i=1}^{3}\int_{\textbf{Y}_{f}}\boldsymbol{W}_{f}^{(i)}(\boldsymbol{y})dy
<\frac{\partial}{\partial{x_{i}}}({t\pi}-{p}^{0})(\boldsymbol{x},t)>,
\,\,\,\,\,\,\,\,\,
\end{multline}
where
\begin{equation}\label{eq5.17}
-\nabla\cdot\big(\mu_{1}\mathbb{D}(y,\boldsymbol{W}_{f}^{(i)})=\boldsymbol{e}^{i},\,\,
\nabla\cdot\boldsymbol{W}_{f}^{(i)}=0,\,\boldsymbol{y}\in\textbf{Y}_{f},\,
\boldsymbol{W}_{f}^{(i)}=0,\,\,\boldsymbol{y}\in\gamma.
\end{equation}
The proof of the existence and uniqueness results for the problem (\ref{eq5.17}) is standard and follows from energy estimates
\begin{equation*}
\boldsymbol{W}_{f}^{(i)}\in\mathbb{W}^{(1,0}_{2}(\textbf{Y}_{f}),\,\,
\int_{{Y}_{f}}(|\boldsymbol{W}_{f}^{(i)}|^{2}+|\mathbb{D}(y,\boldsymbol{W}_{f}^{(i)})|^{2})dy
\leqslant\,{M},\,\,i=1,2,3,
\end{equation*}
which are the result of multiplying equation in (\ref{eq5.17}) by $\boldsymbol{W}_{f}^{(i)}$
summing over $i$ from ${1}$ to ${3}$, integration by parts and use of the embedding theorem (\ref{eq2.33}).

Next, define the constant matrix $(B)^{(w)}(r)$ as
\begin{multline}\label{eq5.18}
(B)^{(w)}(r)=\sum_{i,j=1}^{3}{B}^{(w)}_{i,j}(\boldsymbol{e}^{i}\otimes\boldsymbol{e}^{j}+
\boldsymbol{e}^{j}\otimes\boldsymbol{e}^{i}),\,\,
(B)^{(w)}(r)<\boldsymbol{e}^{(i)},\boldsymbol{e}^{(j)}>={b}^{(w)}_{i,j},
\\
{b}^{(w)}_{i,j}=
\int_{\textbf{Y}_{f}}(\boldsymbol{W}_{f}^{(i)}(\boldsymbol{y})\cdot\boldsymbol{e}^{(j)})dy.
\end{multline}
The matrix $(B)^{(w)}(r)$ is obviously symmetric and strictly positive definite.
In fact, the multiplication equation in (\ref{eq5.18}) by $\boldsymbol{W}_{f}^{(j)}$ results
\begin{equation*}
\mu_{1}\int_{{Y}_{f}}\mathbb{D}(y,\boldsymbol{W}_{f}^{(i)}:
\mathbb{D}(y,\boldsymbol{W}_{f}^{(j)})dy=
\int_{\textbf{Y}_{f}}(\boldsymbol{W}_{f}^{(j)}\cdot\boldsymbol{e}^{(i)}dy.
\end{equation*}
Then the equality
\begin{equation*}
\mu_{1}\int_{\textbf{Y}_{f}}\mathbb{D}(y,\boldsymbol{W}_{f}^{(i)}:
\mathbb{D}(y,\boldsymbol{W}_{f}^{(j)})dy=
\int_{\textbf{Y}_{f}}\mu_{1}\mathbb{D}(y,\boldsymbol{W}_{f}^{(j)}:
\mathbb{D}(y,\boldsymbol{W}_{f}^{(i)})dy
\end{equation*}
implies equality
\begin{equation}\label{eq5.19}
\int_{\textbf{Y}_{f}}(\boldsymbol{W}_{f}^{(i)}\cdot\boldsymbol{e}^{(j)})dy=
\int_{\textbf{Y}_{f}}(\boldsymbol{W}_{f}^{(j)}\cdot\boldsymbol{e}^{(i)})dy,
\end{equation}
that means the symmetry of the matrix $(B)^{(w)}(r)$.

To prove the strict positive definiteness of the matrix $(B)^{(w)}(r)$ we put
\begin{equation*}
\boldsymbol{W}_{f}(\boldsymbol{\xi})=\sum_{i=1}^{3}\boldsymbol{W}_{f}^{(i)}\xi_{i}
\end{equation*}
for any vector $\boldsymbol{\xi}\in\mathbb{R}^{3}$ and any function
$\boldsymbol{\varphi}\in\stackrel{\!\!\circ}{\mathbb{W}}^{1,0}_{2}({Y}_{f})$.

Next we consider the integral identity
\begin{equation}\label{eq5.20}
\int_{\textbf{Y}_{f}}\mathbb{D}\big(y,\boldsymbol{W}_{f}(\boldsymbol{\xi})\big):
\mathbb{D}(y,\boldsymbol{\varphi})dy=
\int_{\textbf{Y}_{f}}(\boldsymbol{\varphi}\cdot\boldsymbol{\xi})dy.
\end{equation}
Then for $\boldsymbol{\varphi}=\boldsymbol{W}_{f}(\boldsymbol{\xi})$ one has
\begin{multline}\label{eq5.21}
\int_{\textbf{Y}_{f}}|\mathbb{D}(y,\boldsymbol{W}_{f}(\boldsymbol{\xi})|^{2}dy=
\int_{\textbf{Y}_{f}}\big(\boldsymbol{\xi}\cdot\boldsymbol{W}_{f}(\boldsymbol{\xi})\big)dy=
\\
\sum_{i,j=1}^{3}{\xi}_{i}{\xi}_{j}\int_{\textbf{Y}_{f}}(\boldsymbol{W}_{f}^{(i)}
\cdot\boldsymbol{e}^{(j)})dy=\sum_{i,j=1}^{3}{b}^{(w)}_{i,j}{\xi}_{i}{\xi}_{j}>0.
\end{multline}
It is evident that the equality $\sum_{i,j=1}^{3}{b}^{(w)}_{i,j}{\xi}_{i}{\xi}_{j}=0$
implies the equalities $\boldsymbol{W}_{f}=0$ and $\mathbb{D}\big(y,\boldsymbol{W}_{f}(\boldsymbol{\xi})\big)=0$ in $\textbf{Y}_{f}$, which is impossible.

Therefore, we can limit ourselves to the case $|\boldsymbol{\xi}|=1$. This fact immediately leads to the inequality
\begin{equation*}
\sum_{i,j=1}^{3}{b}^{(w)}_{i,j}{\xi}_{i}{\xi}_{j}\geqslant\,\alpha_{0}=\mbox{const}>0.
\end{equation*}
\end{proof}

\subsection{Proof of Theorem~\ref{t3.1}: homogenization of the dynamic problem $\mathbb{B}^{\varepsilon}(r)$ for the solid component}

\begin{lemma} \label{l5.3}
Under the conditions of the Theorem~\ref{t2.1} the limiting procedure in the integral identity (\ref{eq2.48}) results in the following dynamic problem $\mathbb{H}(r)$ for displacements and pressure, consisting of the homogenized Lam\'{e}s system
\begin{multline}\label{eq5.22}
\,\,\,\,\,\,\,\,
\nabla\cdot\Big(\lambda_{0}\mathfrak{N}^{(s)}:\big(\mathbb{D}(x,\boldsymbol{w}_{s})+
{c}^{2}_{s}(\nabla\cdot\boldsymbol{w}_{s})\mathbb{I}\big)\Big)=\nabla{p}_{0},\,\,\boldsymbol{x}\in\Omega,
\\
\,\,\,\,\boldsymbol{w}_{s}(\boldsymbol{x},t)=0,\,\,\boldsymbol{x}\in\partial\Omega,\,\,0<t<T,
\end{multline}
where the symmetric strictly positive definite constant tensor $\mathfrak{N}^{(s)}$ in (\ref{eq5.23})
is given by formula (\ref{eq5.27}).
\end{lemma}

\begin{proof}
The limiting procedure in the integral identity (\ref{eq2.47}) for test functions
$\boldsymbol{\varphi}=\boldsymbol{\varphi}(\boldsymbol{x},t)$ gives us
\begin{multline}\label{eq5.23}
\int_{0}^{t_{0}}\int_{\Omega}\Big(\big((\nabla{p}^{0}\cdot\boldsymbol{\varphi})+
\lambda_{0}\mathbb{D}(x,\boldsymbol{w}_{s})+{c}_{s}^{2}(\nabla\cdot\boldsymbol{w}_{s})
\mathbb{I}\big):\mathbb{D}(x,\boldsymbol{\varphi})\Big)dxdt+
\\
\int_{0}^{t_{0}}\int_{\Omega}\int_{\textbf{Y}_{s}}
\Big(\big(\lambda_{0}\mathbb{D}(y,\boldsymbol{W}_{s})+
{c}_{s}^{2}(\nabla\cdot\boldsymbol{W}_{s})\mathbb{I}\big)dy:\mathbb{D}(x,\boldsymbol{\varphi})\Big)dxdt=
I^{(0)}_{s}+I^{(1)}_{s}.
\end{multline}
To calculate the term $I^{(1)}_{s}$, we substitute into identity (\ref{eq2.48}) the term $\displaystyle\,\boldsymbol{\varphi}=\varepsilon\eta(\boldsymbol{x},t)
\boldsymbol{\varphi}_{0}(\frac{\boldsymbol{x}}{\varepsilon})$ and again go to the limit as $\varepsilon\rightarrow{0}$:
\begin{equation}\label{eq5.24}
\int_{0}^{t_{0}}\int_{\Omega}\eta\int_{\textbf{Y}_{s}}
\big(\lambda_{0}\mathbb{D}(y,\boldsymbol{W}_{s})+{c}_{s}^{2}(\nabla\cdot
\boldsymbol{W}_{s})\mathbb{I}\big):\mathbb{D}(y,\boldsymbol{\varphi}_{0})\big)dxdt=0.
\end{equation}
Thus
\begin{equation*}
\nabla_{y}\cdot\big(\lambda_{0}\mathbb{D}(y,\boldsymbol{W}_{s})+{c}_{s}^{2}(\nabla\cdot
\boldsymbol{W}_{s})\mathbb{I}\big)=0,\,\,\,\boldsymbol{y}\in\textbf{Y}_{s},
\end{equation*}
\begin{equation}\label{eq5.25}
\big(\lambda_{0}\mathbb{D}(y,\boldsymbol{W}_{s})+{c}_{s}^{2}
(\nabla\cdot\boldsymbol{W}_{s})\mathbb{I}\big)<\boldsymbol{N}>=0,\,\,\,\boldsymbol{y}\in\gamma.
\end{equation}
To define the tensor $\mathbb{D}(y,\boldsymbol{W}_{s})$ we use the representation
\begin{equation}\label{eq5.26}
\begin{rcases}& \text{$\boldsymbol{W}_{s}(\boldsymbol{y};\boldsymbol{x},t)=
\sum_{i,j=1}^{3}\big(\boldsymbol{W}^{ij}_{s}(\boldsymbol{y})d_{ij}(\boldsymbol{x},t)+
\boldsymbol{W}^{ji}_{s}(\boldsymbol{y})d_{ji}(\boldsymbol{x},t)\big)$},
\\
& \text{$\nabla_{y}\Big(\cdot\big(\mathbb{D}(y,\boldsymbol{W}^{ij}_{s})+
{c}_{s}^{2}\nabla\cdot(\boldsymbol{W}^{ij}_{s})\mathbb{I}+\mathbb{J}^{ij}\big)\Big)=0,\,\,\,
(\boldsymbol{y};(\boldsymbol{x},t)\in{Y}_{s}\times\Omega_{T}$},\,\,
\\
&\text{$\big(\mathbb{D}(y,\boldsymbol{W}^{ij}_{s})+
{c}_{s}^{2}\nabla\cdot(\boldsymbol{W}^{ij}_{s})\mathbb{I}+
\mathbb{J}^{ij}\big)\cdot\boldsymbol{N})=0,\,\,\boldsymbol{y}\in\gamma(r),\,\,\,i,j=1,2,3$}.
\end{rcases}
\end{equation}
Due to the a priori estimate
\begin{equation*}
\int_{0}^{t_{0}}\int_{\Omega}\int_{\textbf{Y}_{s}}\big(|\mathbb{D}(y,\boldsymbol{W}^{ij}_{s}|^{2}+
|\nabla\cdot(\boldsymbol{W}^{ij}_{s})|^{2}\big)dydxdt\leqslant\,M^{2},
\end{equation*}
which is the result of multiplying the equation in (\ref{eq5.26}) by $\boldsymbol{W}^{ij}_{s}$ and integrating by parts and the symmetry of $\boldsymbol{W}^{ij}_{s}$, that follows from the symmetry of $\mathbb{J}^{ij}$, the periodic problem (\ref{eq5.26}) obviously has a unique classical solution, such that $\boldsymbol{W}_{s}^{(ji)}=\boldsymbol{W}_{s}^{(ij)}$, i,j=1,2,3, with a constant $M$ depending on $p^{0}$ (see~\cite{29}).
Finally, we define the tensor $\mathfrak{N}_{s}$ as
\begin{equation}\label{eq5.27}
\mathfrak{N}^{(s)}=\lambda_{0}\sum_{i,j=1}^{3}\int_{\textbf{Y}_{s}}
\mathbb{D}(y,\boldsymbol{W}^{ij}_{s})\otimes\mathbb{D}(y,\boldsymbol{W}^{ij}_{s})dy.
\end{equation}
The symmetry of the tensor $\mathfrak{N}^{(s)}$ is obvious. The strict positive definiteness follows from the equality
\begin{multline*}
\lambda_{0}\sum_{i,j=1}^{3}\int_{\textbf{Y}_{s}}\mathfrak{N}^{(s)}
<\mathbb{D}(y,\boldsymbol{W}^{ij}_{s}),\mathbb{D}(y,\boldsymbol{W}^{ij}_{s})>=
\\
\lambda_{0}\sum_{i,j=1}^{3}\int_{\textbf{Y}_{s}}\mathbb{D}(y,\boldsymbol{W}^{ij}_{s}):
\mathbb{D}(y,\boldsymbol{W}^{ij}_{s})=
\lambda_{0}\sum_{i,j=1}^{3}\int_{\textbf{Y}_{s}}d^{2}(y,\boldsymbol{W}^{ij}_{s})>0.
\end{multline*}
for $\displaystyle\,\sum_{i,j=1}^{3}\int_{\textbf{Y}_{s}}|\boldsymbol{W}^{ij}_{s}|^{2}>0$.

In our case the condition
$\displaystyle\,\sum_{i,j=1}^{3}\int_{\textbf{Y}_{s}}|\boldsymbol{W}^{ij}_{s}|dy=0$
means that $\boldsymbol{W}^{ij}_{s}(\boldsymbol{y})=0$ in $\textbf{Y}_{s}$ and $\mathfrak{N}^{(s)}=0$, which is impossible.

Multiplying equation (\ref{eq5.22}) by $\boldsymbol{w}_{s}$ and integrating by parts gives the result $\boldsymbol{w}_{s}=0$ in $\Omega_{T}$. But this fact contradicts the equation (\ref{eq5.22}), because $|\nabla{p}_{0}|>0$ in $\Omega_{T}$. Thus the tensor $\mathfrak{N}^{(s)}$ is strictly positive definite.
\end{proof}

\subsection{Proof of Theorem~\ref{t3.3}: homogenization of the diffusion problem $\mathbb{B}^{\varepsilon}(r)$}

\begin{lemma} \label{l5.4}
Under the conditions of Theorem~\ref{t2.1} the limiting procedure in the integral identity (\ref{eq2.49}) results in the following homogenized diffusion problem $\mathbb{H}(r)$ for the concentration ${c}\in\mathbb{H}^{2+\alpha,\frac{2+\alpha}{2}}(\overline{\Omega}_{T})$ of the acid, consisting of the differential equation
\begin{equation}\label{eq5.28}
\frac{\partial}{\partial{t}}\big(m(r){c}\big)=
\nabla\cdot(\alpha_{c}\mathbb{B}^{(c)}(r)<\nabla({c}-{c}^{0})>)
\end{equation}
in the domain $\Omega_{T}$ and the boundary and initial conditions
\begin{equation}\label{eq5.29}
{c}(\boldsymbol{x},{t})={c}^{0}(\boldsymbol{x}),\,\,\,\boldsymbol{x}\in{S}^{1}\cup{S}^{2},\,\,
\frac{\partial{c}}{\partial{n}}(\boldsymbol{x},t)=0,\,\,\boldsymbol{x}\in S^{0},\,\,t>0,
\end{equation}
\begin{equation}\label{eq5.30}
{c}(\boldsymbol{x},0)={c}^{0}(\boldsymbol{x}),\,\,\,\boldsymbol{x}\in\Omega.
\end{equation}
The symmetric strictly positive definite constant matrix $\mathbb{B}^{(c)}(r)$ is given by the formula:
\begin{equation}\label{eq5.31}
\mathbb{B}^{(c)}(r)=\frac{1}{2}\sum_{i=1,j}^{3}\int_{\textbf{Y}_{f}}
\big(\nabla_{y}{C}^{i}(\boldsymbol{y})\otimes\nabla_{y}{C}^{j}+
(\nabla_{y}{C}^{j})\otimes\nabla_{y}{C}^{i})\big)dy
\end{equation}
and hold true the following a priori estimates
\begin{equation}\label{eq5.32}
{0}\leqslant{c}(\boldsymbol{x},t)\leqslant{1},\,\,(\boldsymbol{x},t)\in\Omega_{T},\,\,\,
|{c}|^{(2+\alpha)}_{\Omega_{T}}\,\leqslant\,M.
\end{equation}
\begin{proof}
The first estimate in (\ref{eq5.32}) is the maximum principle (consequence~2.2 of Theorem~2.1, chapter~1, \S~2 in~\cite{28}).

To prove the second estimate we consider the integral identity (\ref{eq2.49}) in the form
\begin{equation}\label{eq5.33}
\int_{0}^{t_{0}}\int_{\Omega}\chi^{\varepsilon}\big(-({\widetilde{c}}^{\,\varepsilon}+
\frac{\beta^{\varepsilon}}{\alpha^{\varepsilon}})\,\frac{\partial\xi}{\partial{t}}+
\nabla\xi\cdot(\alpha_{c}\nabla{\widetilde{c}}^{\,\varepsilon}\big)dxdt=0
\end{equation}
for any arbitrary smooth functions $\xi$, vanishing at the boundary $\big(S^{1}\cup S^{2})$,
at $t=0$ and at $t=t_{0}$.

In accordance with Lemma~\ref{l5.1} we get

1) the sequence $\{\widetilde{c}^{\,\varepsilon}\}$ converges weakly and two-scale to the function
${c}\in\mathbb{W}^{1,0}_{2}(\Omega_{T})$;

2) the sequence $\{\nabla\widetilde{c}^{\,\varepsilon}\}$ converges two-scale to the function
${\nabla{c}}+\nabla_{y}{C}$.

Note that weak and two-scale convergence are connected by the relation:
\begin{equation*}
\begin{aligned}
&\text{if } c^{\varepsilon}
\stackrel{2\text{-sc.}}{\longrightarrow}
C(\boldsymbol{x},t;\boldsymbol{y})
\quad \text{(two-scale convergence)}, \\
&\text{then } c^{\varepsilon}(\boldsymbol{x},t)
\rightharpoonup
\int_{\mathbf{Y}_f} C(\boldsymbol{y};\boldsymbol{x},t)\, d\boldsymbol{y}
\quad \text{(weak convergence)}.
\end{aligned}
\end{equation*}
Next, taking into account Lemma~\ref{l5.1} we pass to the limit as $\varepsilon\rightarrow{0}$ in (\ref{eq5.34}) and obtain
\begin{equation}\label{eq5.34}
\int_{0}^{t_{0}}\int_{\Omega}\chi^{\varepsilon}\big((-c+\frac{1}{\theta})\frac{\partial{\xi}}{\partial{t}}+
\alpha_{c}\nabla{\xi}\cdot(\nabla{c}+\int_{{Y}_{f}}\nabla_{y}{C}dy)\big)dxdt=0,
\end{equation}
which results in the diffusion equation for the concentration of acid
\begin{equation}\label{eq5.35}
\frac{\partial}{\partial{t}}\big(m(r){c}\big)=
\alpha_{c}\nabla\cdot\big(\nabla_{x}{c}+\int_{{Y}_{f}}\nabla_{y}{C})dy\big).
\end{equation}
To calculate the integral $\int_{{Y}_{f}}\nabla_{y}{C}dy$ in identity (\ref{eq5.34}), we will again consider the integral identity (\ref{eq5.34}) with test functions $\xi(\boldsymbol{y};\boldsymbol{x},t)=\varepsilon\eta(\boldsymbol{x},t)
\phi(\frac{\boldsymbol{x}}{\varepsilon})$, take the limit as $\varepsilon\rightarrow{0}$ and get the identity
\begin{equation*}
\int_{0}^{t_{0}}\int_{\Omega}\eta(\boldsymbol{x},t)\big(\int_{{Y}_{f}}
\nabla_{y}\phi(\boldsymbol{y})\cdot(\nabla_{x}{c}+\nabla_{y}{C}dy)\big)dxdt=0,
\end{equation*}
which results in the differential equation
\begin{equation}\label{eq5.36}
\nabla_{y}\cdot(\nabla_{x}{c}+\nabla_{y}{C})=0,\,\,\,\boldsymbol{y}\in{Y}_{f}
\end{equation}
for 1-periodic in $\boldsymbol{y}$ functions $C(\boldsymbol{y};(\boldsymbol{x},t))$ and the boundary condition
\begin{equation}\label{eq5.37}
\big((\nabla_{x}{c}+\nabla_{y}{C})\cdot\boldsymbol{N}\big)=0,\,\,\,\boldsymbol{y}\in\gamma,
\end{equation}
where $\boldsymbol{N}$ is the unit normal vector to the boundary $\gamma$.

To solve the problem (\ref{eq5.36}), (\ref{eq5.37}) we use the decomposition
\begin{equation}\label{eq5.38}
\nabla_{y}{C}(\boldsymbol{y};\boldsymbol{x},t)=
\sum_{i=1}^{3}\frac{\partial{C}}{\partial{y}_{i}}(\boldsymbol{y})\boldsymbol{e}_{i},\,\,\,
\nabla_{x}{c}(\boldsymbol{x},t)=\sum_{i=1}^{3}
\frac{\partial{c}}{\partial{x}_{i}}(\boldsymbol{x},t)\boldsymbol{e}_{i}.
\end{equation}
Thus,
\begin{equation}\label{eq5.39}
\sum_{i=1}^{3}\big(\nabla_{y}\cdot(\frac{\partial{C}}{\partial{y}_{i}}+
\frac{\partial{c}}{\partial{x}_{i}})\boldsymbol{e}_{i}\big)=0,\,\,\,\boldsymbol{y}\in\textbf{Y}_{f},
\end{equation}
\begin{equation}\label{eq5.40}
\sum_{i=1}^{3}\big(\frac{\partial{C}}{\partial{y}_{i}}+
\frac{\partial{c}}{\partial{x}_{i}})(\boldsymbol{e}_{i}\cdot\boldsymbol{N})\big),
\,\,\,\boldsymbol{y}\in\gamma.
\end{equation}
The problem (\ref{eq5.39}), (\ref{eq5.40}) obviously has a unique 1-periodic in the variable $\boldsymbol{y}$ classical solution due to the evident a priori estimate
\begin{equation}\label{eq5.41}
\int_{\textbf{Y}_{f}}|\nabla_{y}C|^{2}dy,\leqslant\,M\int_{\Omega_{f}}|\nabla_{x}{c}|^{2}dx
\end{equation}
and estimate (\ref{eq5.2}).
\end{proof}
The matrix $\mathbb{B}^{(c)}(r)$ is obviously symmetric.

To prove the strict positive definiteness of the matrix $\mathbb{B}^{(c)}(r)$ we consider the sum
\begin{multline}\label{eq5.42}
\mathbb{B}^{(c)}(r)<\boldsymbol{\xi}>=\frac{1}{2}\sum_{i=1,j}^{3}
(\nabla_{y}C^{i}(\nabla_{y}C^{j}\cdot\boldsymbol{\xi})+
\nabla_{y}C^{j}(\nabla_{y}C^{i}\cdot\boldsymbol{\xi})),
\\
\mathbb{B}^{(c)}(r)<\boldsymbol{\xi},\boldsymbol{\xi}>=\frac{1}{2}\sum_{i=1,j}^{3}
\Big(\big((\nabla_{y}C^{i}\cdot\boldsymbol{\xi})(\nabla_{y}C^{j}\cdot\boldsymbol{\xi})+
(\nabla_{y}C^{j}\cdot\boldsymbol{\xi})(\nabla_{y}C^{i}\cdot\boldsymbol{\xi})\Big)=
\\
\sum_{i=1,j}^{3}(\nabla_{y}C^{j}\cdot\boldsymbol{\xi})
(\nabla_{y}C^{i}\cdot\boldsymbol{\xi})=
(\sum_{i=1}^{3}\nabla_{y}C^{i}\cdot\boldsymbol{\xi})(\sum_{j=1}^{3}\nabla_{y}C^{j}\cdot\boldsymbol{\xi})=
\\
|\big(\sum_{i=1}^{3}(\nabla_{y}C^{i}\cdot\boldsymbol{\xi})\big)|^{2}>\beta_{0}=\mbox{const}>0
\end{multline}
for an arbitrary vector $\boldsymbol{\xi}$.

In fact, if $\beta_{0}=0$, then for $\displaystyle\,\boldsymbol{\xi}=\sum_{i=1}^{3}\nabla_{y}C^{i}$
we get
\begin{equation*}
|\big(\sum_{i=1}^{3}(\nabla_{y}C^{i})\cdot(\nabla_{y}C^{i})\big)|^{2}=
\sum_{i=1}^{3}|\nabla_{y}C^{i}|^{2}=0,\,\,\mbox{which imply}\,\,\nabla_{y}C^{i}=0,\,\,i=1,2,3,
\end{equation*}
which is impossible.

The smoothness of the function ${c}$ follows from Theorem~5.2, \S~5 in~\cite{28}. In fact, let ${c}^{(1)}=({c}-{c}^{0})$. Then ${c}^{(1)}=0$ at the boundary ${S}^{1}\cup{S}^{2}$ and $(\nabla{c}^{(1)}\cdot\boldsymbol{n})=0$ at the boundary $S^{0}$, where $\boldsymbol{n}$ is the normal vector to the boundary $\partial\Omega$. Therefore we may extend ${c}^{(1)}$ in an even and odd manner across the boundary $S=\partial\Omega$ onto the domain $\Omega^{(1)}$ and get the function ${c}^{(2)}$.

Similarly, we can extend the function ${c}^{(2)}$ across the boundary $\partial\Omega^{(1)}$ onto the domain $\Omega^{(2)}$ and obtain the function ${c}^{(3)}$. It is clear that $\overline{\Omega}\subset\Omega^{(2)}$ and that ${c}^{(3)}$ is a bounded solution of the equation (\ref{eq5.28}) in the domain $\Omega^{(2)}$ and the distance between $\Omega$ and $\Omega^{(2)}$ is positive. Therefore, applying local estimates for parabolic equations (\S~10, Chapter~IV in~\cite{28}), and the maximum principle (the first estimate in (\ref{eq5.32})) we arrive at the second estimate in (\ref{eq5.32}.)
\end{lemma}

\subsection{Homogenization of the boundary condition \eqref{eq2.11}}

\begin{lemma} \label{l5.7}
Let $r\in\mathfrak{M}_{(0,T)}$ and
\begin{equation}\label{eq5.43}
\alpha^{\varepsilon}=\varepsilon{\theta},\,\,  \beta^{\varepsilon}=\varepsilon,
\end{equation}
where $\theta$ is a given positive constant.

Then the homogenization of the boundary condition (\ref{eq2.11}) is given by the formula
\begin{equation}\label{eq5.44}
d_{n}(\boldsymbol{x},t)=\frac{\partial{r}}{\partial{t}}(\boldsymbol{x},t)=
\theta\,c(\boldsymbol{x},t),\,\,\,r(\boldsymbol{x},0)=r_{0}(\boldsymbol{x}).
\end{equation}
For the proof, see Lemma~4.2 in~\cite{18}.
\end{lemma}

\section{Proof of Theorem~\ref{t3.3}: Existence of the Classical Solution to the Problem $\mathbb{H}$}

Let ${c}=\mathbb{F}_{c}(r)$ be the solution to the problem $\mathbb{H}(r)$ and
\begin{equation}\label{eq6.1}
R(\boldsymbol{x},t)=r_{0}(\boldsymbol{x})-
\theta\,\int_{0}^{t}c(\boldsymbol{x},\tau)d\tau\,\equiv \mathbb{F}(r)(\boldsymbol{x},t).
\end{equation}
It is easy to see that the operator $\mathbb{F}(r)$, defined by formula
\eqref{eq6.1} satisfies the Lipschitz condition.
Moreover, for some small time interval $(0,T_{1})$ it is compressive and
displays the set $\mathfrak{M}_{(0,T)}$ into itself (see estimate (\ref{eq3.2})).

In fact, one has:
\begin{multline*}
0\,\leqslant\,\mathbb{F}(r)(\boldsymbol{x},t)\,\leqslant\,T_{1},\,\,\,
|\mathbb{F}(r)|^{(2+\gamma)}_{\Omega_{T}}\,\leqslant\,T_{1}\,M_{0},\\
|\mathbb{F}(r_{1}) - \mathbb{F}(r_{2)}|^{(2+\alpha)}_{\Omega_{T}}
\,\leqslant\,T_{1}\,M_{0}\,|r_{1} - r_{2}|^{(2+\alpha)}_{\Omega_{T}}.
\end{multline*}
That is, on the interval $(0,T_{1})$, where
\begin{equation*}
T_{1}<\frac{M_{0}}{2}.
\end{equation*}
the operator $\mathbb{F}(r)$ is still compressive and maps the set
$\mathfrak{M}_{T_{1}}$ into itself.

Banach's Theorem (Theorem~1, \S~4, Chapter~II, \cite{36}) guarantees
us the existence of the unique fixed point $r_{1}(\boldsymbol{x},t)$
from the set $\mathfrak{M}_{(0,T_{1})}$ and thus establishes the validity of Theorem~\ref{t3.3}
on the time interval $(0,T_{1})$. That is, the existence of the unique solution of the problem $\mathbb{H}(r_{1})$ on the time interval $(0,T_{1})$.

Repeating the procedure, we will get the sequence of the problems $\mathbb{H}(r_{k}),\,\,k=1,2,...$, where
$r_{k}\in\mathfrak{M}_{(0,T)}$, $r_{2}(\boldsymbol{x},T_{1})=r_{1}(\boldsymbol{x},T_{1})$, $r_{3}(\boldsymbol{x},T_{2})=r_{2}(\boldsymbol{x},T_{2})$,...,$r_{k+1}(\boldsymbol{x},T_{k})=
r_{k}(\boldsymbol{x},T_{k})$, $T_{k}<T_{k+1}<T$, $0\leqslant\,r_{k+1}(\boldsymbol{x},t)\leqslant\,
r_{k}(\boldsymbol{x},t)$,\,\,k=1,2,....

There are two possibilities: 1) for some $k^{*}$ $r_{k^{*}}(\boldsymbol{x},t)=r^{*}(\boldsymbol{x},t)$ and $\mathbb{H}(r^{*})=\mathbb{H}$, which solves the problem, and 2) $r_{k}(\boldsymbol{x},t)>0$ for all $k$ and all $t>0$. Obviously, in this case $|r_{k}|^{(0)}_{\Omega_{T}}\rightarrow{0}$ for $k\rightarrow\infty$ and $T_{k}\rightarrow{T}$, which again solves the problem.

\begin{remark}
The similar result has been obtained in \cite{DiazMeirmanov} for compressible components. But in the mentioned publications, all matrices and tensor in   equations (\ref{eq5.32}), (\ref{eq5.32})  (\ref{eq5.32})
\end{remark}
\section*{Acknowledgments}

On behalf of all authors, the corresponding author states that there is no conflict of interest.

\bibliographystyle{siamplain}

\end{document}